\documentclass[onecolumn]{elsart}    

\usepackage{algorithm}
\usepackage{algpseudocode}
\usepackage[table]{xcolor}
\usepackage{graphicx}
\usepackage{epsfig}
\usepackage{amssymb}
\usepackage{amsfonts}

\usepackage{epstopdf}
\usepackage{bbm}
\usepackage{multirow}
\usepackage{tikz}
\usetikzlibrary{patterns,decorations.pathmorphing}
\usepackage{wrapfig}
\usepackage{amsmath}
\usepackage{adjustbox}

\usepackage{caption}
\usepackage{booktabs}
\usepackage{siunitx}
\usepackage{subcaption}

\def\NoNumber#1{{\def\alglinenumber##1{}\State #1}\addtocounter{ALG@line}{-1}}
\renewcommand{\thesubsection}{\Alph{subsection}}
\newcommand{\tN}{[0,N]_{\mathbb Z}}
\newcommand{\Bc}{\mathcal B}
\newcommand{\Nc}{\mathcal N}
\newcommand{\Ns}{\mathbb N}
\newcommand{\nn}{\nonumber}

\newcommand{\tr}{\operatorname{tr}}

\newcommand{\Rs}{\mathbb R}
\newcommand{\Ds}{\mathcal D}

\newcommand{\al}[1]{\begin{align} #1 \end{align}}

\begin{document}
\begin{frontmatter}
\title{A robust approach to sigma point Kalman filtering}

\thanks[footnoteinfo]{This paper was not presented at any IFAC
meeting. Corresponding author: S. Yi.}

\author[Italy]{Shenglun Yi}\ead{shenglun@dei.unipd.it},    
\author[Italy]{Mattia Zorzi}\ead{zorzimat@dei.unipd.it},               

\address[Italy]{Department of Information Engineering, University of Padova, Via Gradenigo 6/B, 35131 Padova, Italy}  



\begin{abstract}
 We propose a robust estimator for nonlinear state-space models and provide a clear interpretation of it as the minimizer of a minimax game. The corresponding maximizer searches for the least favorable model over an ambiguity set whose center is obtained by approximating the nominal model through a sigma-point transformation. Moreover, we develop a Markov Chain Monte Carlo (MCMC) scheme for generating adversarial data from it, thereby allowing the assessment of the resulting uncertainty. \end{abstract}
\end{frontmatter}
\section{Introduction}

Nonlinear state estimation for discrete-time stochastic systems is a classical problem with broad applications in control, signal processing, and autonomous systems. A standard approach is to linearize the model, leading to the extended Kalman filter (EKF) \cite{reif1999stochastic}. However, EKF is effective only for mildly nonlinear systems, since it relies on local linearization based on a first-order Taylor expansion. To achieve a more accurate treatment of nonlinearities while retaining moderate computational complexity, sigma-point Kalman filters, including the unscented Kalman filter (UKF) \cite{julier2000new}, the cubature Kalman filter (CKF) \cite{arasaratnam2009cubature}, and the Gauss--Hermite Kalman filter \cite{pruher2020improved,arasaratnam2007discrete}, have become standard tools for nonlinear estimation \cite{sarkka2023bayesian}.

In many practical applications, however, the estimation problem is characterized not only by nonlinearities, but also by model uncertainty \cite{cheng2026distributionally1,cheng2026distributionally2}. Such uncertainty may arise from imprecisely known model parameters, unmodeled dynamics, or unknown sensor drifts \cite{8025799,9146725}. Under these conditions, standard sigma-point filters may suffer substantial performance degradation, since they are built upon a nominal model. Existing robust sigma-point filtering methods are often specialized to particular types of uncertainty, such as outliers or heavy-tailed non-Gaussian noise \cite{nakabayashi2019nonlinear,zhao2022robust}, rather than to more general perturbations affecting the underlying probabilistic model. A natural framework for handling this broader form of uncertainty is provided by  minimax robust estimation paradigms, which have been extensively studied in the linear setting \cite{ROBUST_STATE_SPACE_LEVY_NIKOUKHAH_2013,STATETAU_2017,kim2020robust,abadeh2018wasserstein,robustleastsquaresestimation,yi2021robust,10654520}. In that framework, the estimator is designed according to the least favorable model selected within a prescribed ambiguity set,   centered at the nominal model, with radius quantifying the level of uncertainty.
The few existing works on robust minimax state estimation for nonlinear models are based on EKF-type approximations \cite{longhini2021learning,11107768,jang2026residual}.  Moreover, the latter mainly focus on the construction of the robust estimator, while they do not characterize the associated least favorable model.  It is worth noting that characterizing the least favorable model is equally important, since it provides the basis for generating adversarial data and assessing whether such a type of uncertainty is realistic.
 
Extending a minimax approach  in  \cite{ROBUST_STATE_SPACE_LEVY_NIKOUKHAH_2013,URKF} to the nonlinear case is highly nontrivial.  Indeed, the challenge is not merely to   counteract uncertainty by enlarging the covariance matrices of the noise processes  in a standard nonlinear algorithm,  but rather to provide a clear interpretation of a robust estimator, as a minimizer of a  minimax game, obtained via a suitable approximation.  Moreover, unlike in the linear case, it is difficult to derive a state-space representation of the least favorable model, i.e. the maximizer, which in turn makes the generation of adversarial data challenging.

In this paper, we consider a minimax game for nonlinear state estimation under model uncertainty in the measurement equation which extends the linear framework recently proposed in \cite{URKF}.  The characterization of the solution to this problem is not tractable. We show that  replacing the center of the ambiguity set with  a sigma-point approximation allows us to characterize a minimizer, thereby yielding a  robust sigma-point-like filter.

Another contribution of the paper concerns the characterization and assessment of the least favorable model, i.e.  the maximizer of the corresponding problem.  The latter does not admit a state-space representation, which  makes direct simulation highly nontrivial. To address this challenge, we propose  a Markov chain Monte Carlo (MCMC) scheme for generating adversarial data from the least favorable model. In particular, we design a proposal state-space model that remains close to the target density while being easy to sample from. We further establish the convergence of the resulting Markov chain.

The  numerical experiments confirm the convergence of the proposed MCMC algorithm and show that the robust sigma-point filter exhibits the expected optimal behavior relative to the standard nonlinear filter under adversarial data. In addition, we consider a more pragmatic setting  with parametric model mismatch and different sensor non-idealities scenarios, where it is not guaranteed that the actual model  belongs to the prescribed ambiguity set. 
Even in this case, the proposed approach compares favorably with the existing state-of-the-art  robust nonlinear filters, such as the  robust EKF \cite{longhini2021learning,11107768}, the maximum correntropy UKF  \cite{zhao2022robust}, and the maximum entropy particle filter \cite{9872130}, in terms of both estimation performance and computational efficiency. Moreover, the simulation results show that the robust EKF, which relies on a simpler approximation of the nonlinear dynamics, yields poor estimation performance. This illustrates that, in practice, it is important not only to guarantee robustness against uncertainty, but also to retain a good approximation of the nonlinear dynamics.  This paper builds upon some preliminary results presented in the conference paper  \cite{yi2026eccspKF}.

The remainder of the paper is organized as follows. Section \ref{sec_2} introduces the problem formulation. Section \ref{sec_3_true} introduces the minimax problem and derives a saddle point solution. Section \ref{sec_lfm} presents the MCMC scheme. Section \ref{sec_8} reports the numerical examples. Finally, Section \ref{sec_conc} concludes the paper.

{\em Notation.} $\tN$ denotes the interval of integers between $0$ and $N$. $z\sim f(z)$ means that the random vector $z$ is distributed according to the probability density $f(z)$; $f(z)= \Nc(\mu ,P)$ means that the probability density $f(z)$ is Gaussian with mean $\mu$ and covariance matrix $P$. Given a sequence of random variables $a_r$, with $r\in\Ns$, $a_r\overset{p}{\longrightarrow} a$ means that the sequence $a_r$ converges in probability to the random variable $a$. Given a symmetric matrix $P$, $P>0$ and $P\geq 0$ mean that $P$ is positive definite and semi-definite, respectively. Moreover, $|P|$ and $\tr(P)$ denote the determinant and the trace of $P$; $I_n$ denotes the identity matrix of dimension $n$; $A^\top$ denotes the transpose of matrix $A$.  Finally, throughout the paper, we omit the arguments of probability density functions, i.e., $\psi_t$ instead of $\psi_t(y_t|x_t)$, when the context is clear in order to keep the notation concise.

\section{Problem Formulation}\label{sec_2}
We consider the nominal  discrete-time nonlinear state space model:
 \begin{equation}\label{nomi_mod}
 \left\{\begin{aligned} x_{t+1}&= f(x_t) + Bv_t \\
 y_t &= h( x_t) +  Dv_t.\\
 \end{aligned}\right.
 \end{equation}
where $x_t\in\mathbb{R}^n$ is the state,  $y_t\in\mathbb{R}^m$ is the  observation, $v_t \in \mathbb{R}^{n+m}$ is white Gaussian noise (WGN) with covariance matrix equal to $I_{n+m}$,  $x_0$ is Gaussian distributed. Matrices $B \in \mathbb{R}^{ n \times (n+m)}$ and $D \in \mathbb{R}^{ m \times (n+m)}$ are full row rank and such that $BD^\top=0$. In plain words, the noise processes $Bv_t$ and $Dv_t$ are assumed independent. 
 Moreover, we assume that  $v_t$ is independent from the initial state $x_0$. Consider  the nominal  state space model (\ref{nomi_mod}). The latter over the finite time interval $\tN$  is characterized by the  nominal density of $Z_N :=\left[\begin{array}{lllllll}{x^\top_{0}} & {\ldots} & {x^\top_{N+1}}, & {y^\top_{0}} & {\ldots}  & {y^\top_{N}}\end{array}\right]^\top $ which is
\begin{equation}\label{nomi1}
p\left(Z_{N}\right)=  p_0\left(x_{0}\right) \prod_{t=0}^{N} {p_{t}}\left(x_{t+1} | x_{t}\right)  \psi_{t}\left(y_{t} | x_{t}\right) ,
\end{equation}
where $p_0(x_0) = \mathcal N(\hat x_0,\tilde P_0)$ denotes the density of $x_0$, and
 \begin{align} \label{pxt1t} p_t(x_{t+1}|x_t) & = \mathcal{N} (f(x_t), ~ BB^\top),\\
\psi_t(y_t|x_t) & = \mathcal{N} (h(x_t), ~ DD^\top). \nn \end{align}
We assume that the actual density of $Z_N$ takes the form
\begin{equation}
\tilde p\left(Z_{N}\right)=p_0\left(x_{0}\right) \prod_{t=0}^{N} {p_{t}}\left(x_{t+1} | x_{t}\right) \tilde \psi_{t}\left(y_{t} | x_{t}\right),\nn
\end{equation}
i.e. the measurement model  -- described by the conditional density $\tilde \psi_t$ -- is different from the nominal one.
We  can measure the mismatch between the actual  and the nominal measurement model at time $t$ by the conditional Kullback–Leibler (KL) divergence
\begin{align*}
\Ds(\tilde{\psi}_{t},& {\psi}_{t}) : =\iint   \tilde \psi_t(y_t | x_t) \tilde p_t(x_t|Y_{t-1})\ln\left(\frac{\tilde \psi_t}{\psi_t}\right)d y_t d x_t
\end{align*}
where $\tilde p_t(x_t|Y_{t-1})$ is the actual a priori density of $x_t$ given $Y_{t-1}:=[\,y_0^\top\,\ldots\, y_{t-1}^\top\,]^\top$. We assume that $\tilde \psi_t$ belongs to the  ambiguity set
\begin{equation}\label{ambi_st}
{\mathcal B}_t:=\left\{\, \tilde{\psi}_{t} \hbox{ s.t. } \Ds(\tilde{\psi}_{t}, {\psi}_{t}) \leq c_{t}\right\}, 
\end{equation}
where $c_t$ is the so-called tolerance. In plain words, we place an upper bound on model uncertainty at each time step, ensuring that uncertainty is not concentrated on specific time steps where the estimator is more vulnerable. This type of ambiguity set naturally arises when the nominal model \eqref{nomi_mod} is estimated from data,  see \cite[Section 2]{URKF} and \cite{hansen2008robustness}. Notice that, $\Bc_t$ is a convex set, moreover it contains conditional densities which are not necessarily Gaussian. 
Our aim is to develop a general approach for   sigma point Kalman filtering, which provides a robust a posteriori state estimator of $x_{t}$ given $Y_t$,  in the case that the actual measurement model is unknown and different from the nominal one  in (\ref{nomi_mod}). 
When the functions 
$f(x)$ and $h(x)$  are linear in 
$x$, the problem has been addressed in \cite{URKF}, where a robust filter can be derived and interpreted within a dynamic minimax game framework. In that setting, the robust estimator can be explicitly characterized as a  Kalman-like filter with an enlarged estimation error covariance, and the associated least-favorable model can also be derived in closed form as a state space model, which allows the systematic generation of data from such a model. However, when the system is nonlinear, although one may heuristically introduce robustness in a given estimation algorithm by enlarging the estimation error covariance, the interpretation of this modification is far from trivial. In particular, it is generally unclear which underlying dynamic game leads to such a solution, and how the corresponding least-favorable model should be constructed in order to generate  adversarial data consistent with that robustness interpretation. In the following we address these  challenges by providing a dynamic game interpretation and characterizing the associated least-favorable model in the nonlinear setting. Throughout the paper, to ease the exposition, we will consider the case in which the nonlinear model (\ref{nomi_mod}) is time-invariant and autonomous, i.e. there is no exogenous input $u_t$. However, the results we present can be extended to the  time-varying case where $f_t(x_t,u_t)$ and $h_t(x_t,u_t)$ depend on an input $u_t$,  see Remark \ref{rem_gen} for more details.

\section{Robust sigma point Kalman filtering}\label{sec_3_true}

Our framework relies on a dynamic minimax game composed by two players: one player, i.e. the state estimator, minimizes the variance of the state estimation error, while the other one, i.e. the nature, selects the least favorable model belonging to a set of possible models about the nominal model (\ref{nomi_mod}). The latter is the prescribed ambiguity set (\ref{ambi_st}).  In this setting, the robust estimation problem characterizing the  estimator of $x_t$ given $Y_t$ is defined as
\begin{equation}
 \label{minimaxrsu2} (\tilde \psi^\star_t, g^\star_t ) =\arg \underset{g_t \in \mathcal{G}_{t}}{\mathrm{min}}\max_{\tilde{\psi}_{t} \in \mathcal{B}_{t}} J_t(\tilde {\psi}_t,g_t)
\end{equation}
where
{\small \begin{equation}\label{JtU} \begin{aligned}
J_t(\tilde {\psi}_t,g_t)=&\frac{1}{2} \iint \left\|  x_{t}-g_{t}\left(y_{t}\right) \right\|^{2}  \tilde{\psi}_{t}(y_t | x_t)   \tilde{ p}_{t}\left(x_{t} | Y_{t-1}\right) d y_{t} d x_{t}
\end{aligned}
\end{equation}}
is the variance of the estimation error under the actual model; $\mathcal G_t$ is the convex set of estimators with finite second order moments which are square integrable  with respect to densities $\tilde\psi_t \tilde { p}_{t}(x_{t} | Y_{t-1}) $ and $\tilde \psi_t$ satisfies the condition
\begin{equation}
\iint  \tilde{\psi}_{t}(y_t | x_t)\tilde { p}_{t}\left(x_{t} | Y_{t-1}\right)  d y_{t} d x_{t}=1. \label{Iuukf}
\end{equation}
Our aim is to characterize a saddle point solution to \eqref{minimaxrsu2}.  Solutions of this type of problems are typically non-unique, see \cite{ROBUST_STATE_SPACE_LEVY_NIKOUKHAH_2013,hansen2008robustness} for instance. In view of the Von Neumann's minimax theorem \cite{aubin2006applied},  there exists a saddle point $(\tilde \psi^{\star}_t,~ g^{\star}_t)$ such that
\begin{equation}
                    J_t(\tilde \psi_t,~ g^{\star}_t) \leq J_t(\tilde \psi^{\star}_t,~ g^{\star}_t) \leq J_t(\tilde \psi^{\star}_t,~ g_t),\nn
                  \end{equation}
since the corresponding sets $\mathcal{B}_t$ and $\mathcal{G}_t$ are convex and compact. The next result characterizes the structure of the maximizer of (\ref{minimaxrsu2}), i.e. the least favorable model.
\begin{prop} \label{l_lfm2}
For a fixed estimator $g_t \in \mathcal G_t$, if there exists a density
  \begin{equation}\label{psi_opt}
 \tilde{\psi}_t (y_t|x_t) =\frac{1}{\kappa_t } \exp \left(\frac{ \theta_{t|t}}{2} \left\|x_{t}- g_{t}(y_t)\right\|^2\right)  \psi_t(y_t|x_t),
\end{equation}
where $\theta_{t|t} > 0$ is the unique solution to $\mathcal{D}(\tilde \psi_t,  \psi_t) = c_t$, and the normalizing constant is given by \begin{equation} \label{MT1}
\begin{aligned} 
\kappa_t  = \iint \exp & \left(\frac{\theta_{t|t}}{2} \| x_{t}- g_{t}(y_t) \|^2\right)  \psi_t\tilde{p}_t\left(x_t|Y_{t-1}\right)    d y_t d x_t,
\end{aligned}
\end{equation}
then  $\tilde \psi_t\in \mathcal B_t$ and it is the maximizer of (\ref{JtU}).
\end{prop}
\begin{pf} The claim can be proved using the method of Lagrange multipliers, similarly to \cite[Lemma 2]{URKF}. Although the result in \cite{URKF} was established in the context of a linear state-space model, the derivation does not rely on the linearity assumption. However, unlike the linear case, the existence of a maximizer is not always guaranteed; as we will see, it depends on the behavior of the function $h$ in \eqref{nomi_mod}.  \qed
\end{pf}

Let 
\al{\bar p_t(w_t|Y_{t-1})&= \psi_t(y_t|x_t)\tilde p_t(x_t|Y_{t-1})\nn\\
\tilde p_t(w_t|Y_{t-1})&= \tilde \psi_t(y_t|x_t)\tilde p_t(x_t|Y_{t-1})\nn
} 
denote the pseudo-nominal and the actual conditional density of  $w_t:=[\,x_t^\top\;
 y_t^\top\,]^\top$ given $Y_{t-1}$, respectively. The idea is to characterize a saddle point solution by mapping Problem \eqref{minimaxrsu2} into a  minimax game whose players are  $g_t$ and $\tilde p_t(w_t|Y_{t-1})$ 
 \al{\label{J_eqpb_minimax}    \underset{g_t \in \mathcal{G}_{t}}{\mathrm{min}}\max_{\tilde{p}_{t} \in \tilde{\mathcal{B}}_{t}}\bar J_t(\tilde {p}_t,g_t)}
 where 
 \al{\bar J_t(\tilde p_t,g_t):=&\frac{1}{2} \iint \left\|  x_{t}-g_{t}\left(y_{t}\right) \right\|^{2}     \tilde{ p}_{t}\left(w_{t} | Y_{t-1}\right) d w_{t};  }
  the ambiguity set $\tilde \Bc_t$ is obtained using the correspondence between  $\tilde \psi_t(y_t|x_t)$ and  $\tilde p_t(w_t|Y_{t-1})$ given by the optimality condition \eqref{psi_opt}:  
 \begin{equation} \label{ball_uukf2}
\tilde{\mathcal B}_t:=\left\{\, \tilde{p}_{t}(w_t|Y_{t-1}) \hbox{ s.t. } \Ds_{KL}(\tilde{p}_{t},\bar {p}_{t}) \leq c_{t}\right\}
\end{equation}
where 
\al{\Ds_{KL}(\tilde{p}_{t},\bar {p}_{t})=\int\ln\left(\frac{\tilde p_t}{\bar p_t}\right)\tilde p_t dw_t.   \nn}
Unfortunately, in this case it is not possible to characterize a saddle point solution to (\ref{J_eqpb_minimax}) because the pseudo-nominal density  $\bar p_t(w_t|Y_{t-1})$ is not Gaussian. Thus, we construct an approximation of $\bar p_t(w_t|Y_{t-1})$ using the same approximation mechanism exploited in   the  traditional sigma point Kalman filter, \cite{sarkka2023bayesian}. First, we introduce the definition of sigma points, which is used to compute an approximation of the expectation of a nonlinear transformation of a Gaussian random vector. 
\begin{defn} Given a random vector $ x\sim \mathcal N(\hat x,P )$, we denote its sigma points as $ \mathcal{X}^i = \sigma_{i}(\hat x,P)$, with $i = 1, \ldots, p,$ and the corresponding weights for the mean and the covariance matrix
 are denoted by $W^i_m$ and $W^i_c \geq 0$.\end{defn}

There are many ways to define the sigma points and the weights (see \cite[Sec. 5-7]{sarkka2023bayesian}). Different choices result in various nonlinear Kalman filters, among which the most popular are:
\begin{itemize}
  \item \textbf{Unscented Kalman filter (UKF).} The sigma points are obtained through the unscented transformation:
{\small
\begin{equation}\nn
\begin{aligned}
\sigma_{i}(\hat x,P)=\left\{\begin{array}{ll}
 \hat x  + \sqrt{\lambda +n} (\sqrt{{P}} )_i, &  \hbox{if } 1\leq i\leq n  \\
\hat x  - \sqrt{\lambda +n} (\sqrt{{P}} )_{i-n}, &  \hbox{if } n+1\leq i\leq 2n\\
\hat x, & \hbox{if } i=2n + 1    \end{array}\right.
 \end{aligned}
\end{equation}}where $(\sqrt{ P })_i\in\Rs^n$
is the $i$-th column of $\sqrt{ P }$ which is a square root matrix of $P$.
The corresponding weights are
{\small
\begin{align}\label{W_ukf}
  W_m^{i}&=\lambda / (n+\lambda), ~ W^i_c = W^i_m +1- a^2 +b, \quad  \hbox{if } i=2n + 1 \nn \\
   W_m^i&=  W_c^i = 1 / (2(n+\lambda)), \quad  \hbox{if } 1\leq i \leq 2n
\end{align}}where $\lambda = a^2 (\kappa +n) - n$, and the parameters $a$, $b$ and $\kappa$ can be chosen as suggested in \cite{wan2001unscented}.
  \item  \textbf{Cubature Kalman filter (CKF).} The sigma points are generated using the spherical cubature  transformation:
\begin{equation}\nn
\begin{aligned}
\sigma_{i}(\hat x,P)=\left\{\begin{array}{ll}
 \hat x  + (\sqrt{{n P}} )_i, &  \hbox{if } 1\leq i\leq n  \\
\hat x  -  (\sqrt{{n P}} )_{i-n}, &  \hbox{if } n+1\leq i\leq 2n  \end{array}\right.
 \end{aligned}
\end{equation}
and the corresponding weights  are $W^i_c = W_m^i  = {1}/{2n}$.
  \item  \textbf{Gauss-Hermite Kalman filter.}  The sigma points are generated by the Gauss-Hermite moment transformation
\begin{equation}\nn
\begin{aligned}
\sigma_{i}(\hat x,P)=\hat x +\sqrt{{P}} \lambda^i,  \quad i= 1\ldots q^n
 \end{aligned}
\end{equation}
where $ \lambda^i \in \Rs^{n}$ is the $i-$th vector of the set formed by the $n-$dimensional Cartesian products of the roots,  say $\nu_{k}$ with $k=1, \ldots, q$, of the Hermite polynomial of order $q$, say $H_q(\nu)$. The
corresponding weights $W_m^i = W_c^i$ are formed as the products of the $n$ terms
$$\begin{aligned}  \frac{q!}{q^2\left(H_{q-1}\left(\lambda_{k}^i\right)\right)^2}
\end{aligned}$$
where $\lambda_{k}^i$, with $k = 1, \ldots, n$, denotes the $k$-th element of $\lambda^i$.
\end{itemize}
 The approximation $\bar  {p}_{t}^a(w_{t} | Y_{t-1})=\mathcal N(m_t,K_t)$ of $\bar  {p}_{t}(w_{t} | Y_{t-1})$ in standard sigma point Kalman filters is obtained from the sigma points corresponding to the approximation $\tilde  {p}_{t}(x_{t} | Y_{t-1})\simeq\mathcal N(\hat x_t,P_t)$:
\al{ \label{defmK}m_t&= \left[\begin{array}{c}  \hat{x}_t \\ m_{y_t}\end{array}\right],\;   K_t=  \left[\begin{array}{cc} P_t & K_{x_t y_t} \\  K_{y_t x_t}  & K_{y_t}\end{array}\right], }
{\small  \begin{equation}\label{mmy}
\begin{aligned}
m_{y_t} &= \sum_{i=1}^{p} W_m^i h(\mathcal{X}^i_t)\\
K_{y_t} &= \sum_{i=1}^{p} W_c^i(h(\mathcal{X}^i_t) - m_{y_t})(h(\mathcal{X}^i_t) - m_{y_t})^\top + D D^\top \\
K_{x_t y_t} &= \sum_{i=1}^{p} W_c^i(\mathcal{X}^i_t - \hat x_t)(h(\mathcal{X}^i_t) - m_{y_t})^\top\\
\mathcal{X}^i_t&=\sigma_i(\hat x_t, P_t), \quad i=1\ldots p.
\end{aligned}\end{equation}}Applying such approximation in \eqref{J_eqpb_minimax} we obtain the  approximate problem  \begin{equation}
 \label{minimaxa_up}    \hat x_{t|t} = \underset{g_t \in \mathcal{G}_{t}}{\mathrm{argmin}}\max_{\tilde  {p}_{t}  \in {\tilde{\mathcal{B}}}^a_{t}} {\bar J_t}(\tilde p_{t}  ,g_t)
\end{equation}
where $\tilde{\mathcal B}^a_t$ is obtained from (\ref{ball_uukf2}) replacing  $\bar{p}_{t}(w_t|Y_{t-1})$ with $\bar{p}_{t}^a(w_t|Y_{t-1})$. In plain words, the approximation regards the center of (the ball describing) the ambiguity set.  Before characterizing a saddle point solution to Problem \eqref{minimaxa_up}, we must specify how $\hat x_t$ and $P_t$, i.e. the mean and covariance matrix of the a priori conditional density at time $t-1$ used in the prediction stage, are obtained. These quantities are computed from: i) the mean and covariance matrix, $\hat{x}_{t-1|t-1}$ and $\tilde{P}_{t-1|t-1}$, of the posterior density $\tilde{p}_{t-1}(x_{t-1}|Y_{t-1})$ at time $t-1$, which is assumed to be reasonably approximated as a Gaussian density; and ii) the state transition model $p_t(x_{t+1}| x_t)$. Since $p_t(x_{t+1}| x_t)$ is not affected by uncertainty, i.e. it represents both the nominal and the actual state process model, we can employ the same approximation used in the standard sigma-point Kalman filter to compute $x_t$ and $P_t$, \cite{sarkka2023bayesian}:
\al{  {\mathcal{X}}_{t-1|t-1}^i&=\sigma_i(\hat x_{t-1|t-1},\tilde P_{t-1|t-1}),\quad i=1\ldots p\nn\\
 \hat{\mathcal{X}}_{t}^i&=f(  {\mathcal{X}}_{t-1|t-1}^i),\quad i=1\ldots p\nn\\
\hat x_t&=\sum_{i=1}^{p} W_m^i\hat{\mathcal X}_t^i,\;  \nn \\
\label{defPP}{ P}_{t}&=\sum_{i=1}^{p} W_c^i(\hat{\mathcal{X}}_{t}^i-\hat  x_{t})(\hat{\mathcal{X}}_{t}^i-
    \hat {x}_{t})^\top + B B^\top.}
We are ready to characterize a saddle point solution to Problem~\eqref{minimaxa_up}.
\begin{thm}
\label{th2}
Let $(\hat x_{t}, P_{t})$ be the prediction pair at time $t-1$ computed as in \eqref{defPP}. The robust estimator solution to (\ref{minimaxa_up}) is
\begin{equation}\label{rob_UPP}
\hat x_{t|t} =\hat x_{t}+{K}_{x_t y_t} {K}_{y_t}^{-1}\left({y}_t-{m}_{y_t}\right).
\end{equation}
The nominal and the least favorable covariance matrix of the estimation error $x_t-\hat x_{t|t}$ are
\begin{align*}
{P}_{t|t} &= {P}_{t}-{K}_{x_t,y_t} {K}_{y_t}^{-1} {K}_{x_t,y_t}^{\top},\\
\tilde P_{t|t}& =   ( P_{t|t}^{-1} -  \theta_{t|t} I) ^{-1}
\end{align*}
where the risk sensitivity parameter $ \theta_{t|t}>0$ is the unique solution to $\gamma( P_{t|t}, \theta_{t|t}) = c_t$  with
\begin{equation} \label{theta}
    \gamma(  P, \theta) := \frac{1}{2}\left(\log|I - \theta  P| + \tr\left((I-\theta P)^{-1} - I \right)\right).
\end{equation}
Moreover, the approximate least favorable density, i.e. 
the maximizer of \eqref{minimaxa_up}, is such that $\Ds_{KL}(\tilde{p}_{t},\bar {p}^a_{t})=c_t$ and the corresponding
posterior density of $x_t$ given $Y_t$  is
\begin{equation}\label{approx2}\tilde p_{t}^a\left({x}_{t} | Y_{t}\right)  =\mathcal{N}( \hat x_{t|t},\tilde  P_{t|t}).\end{equation}
\end{thm}
\begin{pf} First, we show that $K_t$ defined in \eqref{defmK} is positive definite. Let \al{&\mathbf X_t:=[\, \mathcal X_t^1-\hat x_t \,\ldots \, \mathcal X_t^q-\hat x_t]\nn\\
&\mathbf Y_t:=[\, h(\mathcal X_t^1)-m_{y_t} \,\ldots \, h(\mathcal X_t^q)-m_{y_t}]\nn.}
By \eqref{mmy}, we have that $\mathbf X_t=\sqrt{P_t}\Lambda$ where $\Lambda\in\mathbb R^{n\times p}$ and its definition depends on the type of transformation. More precisely: $\Lambda=\sqrt{\lambda+n}[\,I_n\; -I_n\;  0\,]\in\mathbb R^{n\times 2n+1}$ for the unscented transformation; $\Lambda=\sqrt{n}[\,I_n\; -I_n\,]\in\mathbb R^{n\times 2n}$ for the spherical cubature transformation; $\Lambda=[\,\lambda^1\ldots  \lambda^{q^n}\,]\in\mathbb R^{n\times q^n}$ for the Gauss-Hermite moment transformation. It is not difficult to see that $\Lambda W\Lambda^\top=I_n$ where $W$ is the diagonal matrix with entries in the main diagonal $W_c^1\ldots W_c^p$ which denote the corresponding weights. Accordingly, we have $P_t=\mathbf X_t W\mathbf X_t^\top$, $K_{y_t}=\mathbf Y_t W\mathbf Y_t+DD^\top$ and $K_{x_t y_t}= \mathbf X_t W\mathbf Y_t^\top$. Hence,
\al{ K_t= \left[\begin{array}{cc} \mathbf X_t  W \mathbf X_t^\top & \mathbf X_t  W \mathbf Y_t ^\top \\  \mathbf Y_t  W \mathbf X_t ^\top  & \mathbf Y_t  W \mathbf Y_t ^\top +DD^\top\end{array}\right]>0 \nn}
because, by \eqref{defPP}, $\mathbf X_t  W \mathbf X_t^\top\geq BB^\top>0$ and the Schur complement of block $\mathbf X_t  W \mathbf X_t^\top$ of $K_t$ is
\al{\mathbf Y_t&  W \mathbf Y_t ^\top +DD^\top- \mathbf Y_t  W \mathbf X_t^\top (\mathbf X_t  W \mathbf X_t^\top)^{-1}\mathbf X_t  W \mathbf Y_t^\top\nn\\
&\geq  DD^\top>0.\nn}
Since $\bar p_t^a(w_t|Y_{t-1})=\Nc(m_t,K_t)$, with $K_t>0$, we can apply \cite[Theorem 1]{robustleastsquaresestimation}: it follows that the minimizer of (\ref{minimaxa_up}) is (\ref{rob_UPP}) and the least favorable density is $\tilde p_t(w_t|Y_{t-1})=\Nc(m_t,\tilde K_t)$ where
\al{\tilde K_t:= \left[\begin{array}{cc}
\tilde K_{x_{t}} &  K_{x_{t}y_t}  \\
\bar K_{y_tx_{t}}  &  K_{y_{t}}
\end{array}\right]>0, \nn}
$\tilde K_{x_{t}} = \tilde P_{t|t}-K_{x_{t}y_t}K_{y_{t}}^{-1}K_{x_{t}y_t}^\top$  and such that $\Ds_{KL}(\tilde{p}_{t},\bar {p}^a_{t})=c_t$.  Thus, $$\tilde p _t^a(x_t|Y_t)=\frac{\tilde p _t(w_t|Y_t)}{\int \tilde p _t(w_t|Y_t) d x_t}=\mathcal{N}( \hat x_{t|t},\tilde  P_{t|t}).$$ \qed
 \end{pf}
The theorem above not only characterizes the robust estimator, but also shows that the least favorable posterior density, i.e. the one corresponding to \eqref{J_eqpb_minimax}, can be approximated by the Gaussian density  $\tilde{p}_t^a(x_t| Y_t)$ in (\ref{approx2}). The latter allows us to conclude that the sigma-point approximation is valid over the entire time horizon, and this can be shown by an inductive argument: if the posterior $\tilde p_{t-1}(x_{t-1}|Y_{t-1})$ at time $t-1$ can be approximated as Gaussian, then the prior density $\tilde{p}_t(x_t| Y_{t-1})$ can be approximated by means of the  sigma point transformation, and the theorem ensures that the updated posterior at time $t$ remains approximately Gaussian. Consequently, at the next time step, the approximation of $\tilde{p}_{t+1}(x_{t+1} \mid Y_t)$ can again be constructed using sigma points.
\begin{algorithm}[t]
    \caption{Robust sigma point KF at time $t$}\label{RUKF2}
    \hspace*{\algorithmicindent} \textbf{Input} $\hat{x}_{t}$, $ P_{t}$, $c_t$, $y_t$\\
    \hspace*{\algorithmicindent} \textbf{Output}  $\hat{x}_{t+1}$, $ P_{t+1}, \theta_{t|t}$
    \begin{algorithmic}[1]
    \State $\mathcal{X}^i_t=\sigma_i(\hat x_t,  P_t), \quad i=1\ldots p$
    \State  $m_{y_t}=\sum_{i=1}^{p} W_m^i h(\mathcal{X}^i_t)$
    \vspace{0.08cm}
    \State {\small$ {K}_{y_t}=\sum_{i=1}^{p} W_c^i(h(\mathcal{X}^i_t) - m_{y_t})(h(\mathcal{X}^i_t) - m_{y_t})^\top + D D^\top$}
    \vspace{0.08cm}
    \State ${K}_{x_t y_t}=\sum_{i=1}^{p} W_c^i(\mathcal{X}^i_t-\hat x_t)(h(\mathcal{X}^i_t) - m_{y_t})^\top$
        \vspace{0.05cm}
    \State $L_t  = {K}_{x_t y_t} {K}^{-1}_{y_t}$
    \State $\hat x_{t|t}=\hat x_{t}+L_t\left({y}_t-m_{y_t}\right)$
    \State ${P}_{t|t}  ={ P}_{t}-L_t {K}_{y_t} L_t^{\top} $
    \State Find $\theta_{t|t}$ s.t. $\gamma(P_{t|t}, \theta_{t|t}) = c_t$ \label{AL_theta}
    \State $\tilde P_{t|t} = (P_{t|t}^{-1} -  \theta_{t|t} I) ^{-1}$ \label{modU}
\State $\mathcal{X}^i_{t|t}=\sigma_i(\hat x_{t|t},{\tilde {P}_{t|t}}),\quad i=1\ldots p $
        \State   $\hat{\mathcal{X}}^i_{t+1}=f(\mathcal{X}^i_{t|t}),\quad i=1\ldots p $
   \State $\hat {x}_{t+1}=\sum_{i=1}^{p} W_m^i \hat{\mathcal{X}}_{t+1}^i$
    \vspace{0.08cm}
    \State {${ P}_{t+1}=\sum_{i=1}^{p} W_c^i(\hat{\mathcal{X}}_{t+1}^i-\hat  x_{t+1})(\hat{\mathcal{X}}_{t+1}^i-
    \hat {x}_{t+1})^\top$ $\hspace*{1.5cm} + B B^\top$}
    \end{algorithmic}
\end{algorithm} The resulting robust sigma point Kalman filter is outlined in Algorithm \ref{RUKF2}.   Note that, the computation of $\theta_{t|t}$ in Step \ref{AL_theta} can be efficiently computed through  a bisection method, see \cite{zenere2018coupling}. Finally, in the limit case $c_t = 0$, i.e.  the absence of uncertainty in the measurement model at time $t$, we have that $\theta_{t|t}=0$ and $\tilde P_{t|t}=P_{t|t}$,  thus  the robust sigma point Kalman filter coincides with the standard  sigma point Kalman filter. 
\begin{rem}\label{rem_gen}
Consider the case in which the nominal model is time-varying and depends on the input $u_t$, i.e., $f(x_t)$ and $h(x_t)$ in \eqref{nomi_mod} are replaced by $f_t(x_t,u_t)$ and $h_t(x_t,u_t)$, respectively. Assume that $u_t$ is a deterministic signal or a function of the past outputs (i.e., we are considering a model in feedback configuration). At time $t-1$, given the observations $Y_{t-1}$, the inputs $u_t$ and $u_{t-1}$ are known, and therefore $\bar{p}_{t}^a(w_{t} | Y_{t-1})$ is obtained from \eqref{mmy} and \eqref{defPP} by replacing $h(\mathcal{X}_t^i)$ and $f(\mathcal{X}^i_{t-1|t-1})$ with $h_t(\mathcal{X}_t^i,u_t)$ and $f_t(\mathcal{X}^i_{t-1|t-1},u_{t-1})$, respectively. Accordingly, the derivation above still holds, and in Algorithm~\ref{RUKF2} it is sufficient to replace $h(\mathcal{X}_t^i)$ and $f(\mathcal{X}^i_{t|t})$ with $h_t(\mathcal{X}_t^i,u_t)$ and $f_t(\mathcal{X}^i_{t|t},u_t)$, respectively.
\end{rem}

In view of   Theorem \ref{th2}, the approximate least favorable density  takes the following form:
\al{\label{psi0_t}
&\tilde{\psi}_{t}^{0}(y_t|x_t)=\frac{1}{M_{t}} \exp \left(\frac{ \theta_{t|t}}{2} \left\|  x_{t}- \hat x_{t|t}\right\|^{2}\right) \psi_{t}(y_t|x_t)}
 where $\hat x_{t|t}$ and $\theta_{t|t}$ are the state filter and the risk sensitivity parameter obtained by the robust sigma point Kalman filter; \al{
\label{M_T3}
& M_{t}=\iint  \exp \left( \frac{\theta_{t|t}}{2}  \left\|  x_{t}- \hat x_{t|t} \right\|^{2}\right) \psi_{t}  \tilde p^a_{t}\left(x_{t} | Y_{t-1}\right)   d y_{t} d x_{t}}
where 
$$ \tilde p^a_t(x_t|Y_{t-1})= \int \bar p^a_t(w_t|Y_{t-1}) d y_t =\Nc(\hat x_t,P_t)$$
is the marginal obtained from the approximation of the pseudo-nominal density  $\bar p_t(w_t|Y_{t-1})$. It is worth noting that the optimality of $\tilde \psi^0_t$ is guaranteed by Proposition~\ref{l_lfm2} since by direct substitution we have 
\al{\mathcal{D}( \tilde \psi_t^0,  \psi_t) =\Ds_{KL}(\tilde{p}_{t},\bar {p}^a_{t})=c_t \nn}
where $\tilde p_t$ denotes the maximizer of \eqref{minimaxa_up}.  We now investigate the conditions under which the  approximate maximizer is well defined. \begin{thm} 
 Let 
     \al{\label{h_cond}
   \delta_h:= \lim_{\|x\|\to\infty} \frac{\|h(x)\|^2}{\|x\|^2}.
  }
  Assume that $h(x)$ is such that  $\delta_h$ is finite.
Taking the tolerance $c_{t}$ sufficiently small, then the normalization constant $M_t$ in (\ref{M_T3}) is finite, i.e. $\tilde{\psi}_{t}^{0}$ is well defined.
\label{propos2}
\end{thm}
\begin{pf} 
By \eqref{M_T3} we have that
\al{\label{N11}M_t=  \int N_t(x_t) \tilde p_t^a(x_t|Y_{t-1})dx_t }with  
\al{ N_t(x_t) &:=  \int    \exp\left(\frac{ \theta_{t|t}}{2} \| x_{t} - \hat x_{t|t}(y_t)\|^2\right) \psi_t( y_t| x_t)  dy_t, \label{Nxt2}\\ 
 \hat x_{t|t}(y_t)&:=  \hat x_t + L_t(y_t-m_{y_t})\nn}
is the robust estimator in  (\ref{rob_UPP}) with $\theta_{t|t}$, $\hat x_{t}$, $L_t$, $m_{y_t}$ 
computed through Algorithm \ref{RUKF2} using $Y_N$.   By integrating with respect to $y_t$, we obtain 
\begin{equation}\label{hat_M_2}
N_{t}(x_t) = \frac{1}{\sqrt{|S_t|| R|}} 
\exp\Bigg(\frac{1}{2} \big(s_t(x_t)^\top S_t^{-1} s_t(x_t) + l_t(x_t) \big)\Bigg)
\end{equation}
where \begin{equation} \label{slSR}\begin{aligned}
l_t(x_t) &= \theta_{t|t} \| x_t - \hat x_t + L_t m_{y_t} \|^2 - \|  h(x_t)\|^2_{R^{-1}},\\
s_t (x_t)&= \theta_{t|t} L_t^\top (\hat x_t - x_t - L_t m_{y_t}) + R^{-1} h(x_t), \\
R &= D D^\top; ~~ S_t = R^{-1} - \theta_{t|t} L_t^\top L_t. 
\end{aligned}
\end{equation}
Let $\mu(x_t) = [\,(x_t + L_t m_{y_t} - \hat x_t)^\top\; ~ h(x_t)^\top\,]^\top$, then
it is not difficult to see that \al{\label{N22}
N_t(x_t) = \frac{1}{\sqrt{|S_t ||R|}} 
\exp\left(\frac{1}{2} \mu(x_t)^\top \Theta_t \, \mu(x_t)\right),
}
where $\Theta_t$ does not depend on $x_t$ and is defined as
\begin{equation}\label{THETA}
\Theta_t =
\begin{bmatrix}
\theta_{t|t} I + \theta_{t|t}^2 L_t S_t^{-1} L_t^\top & - \theta_{t|t} L_t S_t^{-1} R^{-1} \\
- \theta_{t|t} R_t^{-1} S_t^{-1} L_t^\top & R^{-1} S_t^{-1} R^{-1} - R^{-1}
\end{bmatrix}.
\end{equation} By \cite[Proposition 4.3]{zorzi2017convergence} we know that $\theta\mapsto \gamma (P,\theta)$, with $P>0$, is monotone increasing over $\Rs_+$. Accordingly, the mapping $c_t\mapsto \theta_{t|t}$ characterized  through the equality $c_t=\gamma(P_{t|t},\theta_{t|t})$, is monotone increasing. Therefore, it is possible to obtain $\theta_{t|t}$ arbitrarily small taking $c_t$ sufficiently small.  
The first-order Taylor expansion of $S_t^{-1}$ around $\theta_{t|t} = 0$ is 
\al{\label{defSi}
S_t^{-1} = R + \theta_{t|t} R L_t^\top L_t R +  \mathcal{O}(\theta_{t|t}),
}
where $\mathcal{O}(\theta_{t|t})$ are the infinitesimal  terms of higher order than $\theta_{t|t}$.
Substituting this expansion into \eqref{THETA} and collecting the terms up to first order in $\theta_{t|t}$ yields
\[
\Theta_t =
\theta_{t|t}
\begin{bmatrix}
I & -L_t \\
- L_t^\top & L_t^\top L_t
\end{bmatrix}
+
\mathcal{O}(\theta_{t|t}^2).
\]
Thus,  $\Theta_t \to 0$ as $c_t\to 0$ (and thus $\theta_{t|t} \to 0$). Moreover, $
\Theta_t$ is positive semi-definite for $c_t$ taken sufficiently small.  
We proceed to prove that $M_t$ is finite. By \eqref{N11} and \eqref{N22} we have
\al{\begin{aligned} M_t 
&=  \frac{1}{\sqrt{(2 \pi)^n |P_t|| S_t|| R|}} \times\\
 & \hspace {0.2cm}\int \exp(\frac{1}{2}\mu(x_t)^\top \Theta_t \mu(x_t)  - \frac{1}{2} \|x_t - \hat x_t \|^2_{P^{-1}_t})dx_t .
\end{aligned} \label{inttt}}
Let $\lambda_{max,\Theta_t}$ denote the maximum eigenvalue of $\Theta_t$, and $\lambda_{max,P_t}$ the maximum eigenvalue of $P_t$. Notice that, $\lambda_{max,\Theta_t}\geq0$ because $\Theta_t\geq 0$ for $c_t$ sufficiently small, and $\lambda_{max,P_t}>0$ because $P_t\geq BB^\top >0$. Then, we have 
\al{&\lim_{\|x_t\|\rightarrow \infty} \frac{|\mu(x_t)^\top  \Theta_t \mu(x_t) |}{ \|x_t - \hat x_t \|^2_{P^{-1}_t}}=\lim_{\|x_t\|\rightarrow \infty}\frac{\mu(x_t)^\top  \Theta_t \mu(x_t) }{ \|x_t  \|^2_{P^{-1}_t}} \nn\\
& \leq  \lambda_{max, \Theta_t} \lambda_{max, P_t}\lim_{\|x_t\|\rightarrow \infty}\frac{\|\mu(x_t)\|^2 }{ \|x_t  \|^2}\nn\\
&=\lambda_{max, \Theta_t}\lambda_{max, P_t}\lim_{\|x_t\|\rightarrow \infty}\frac{\|x_t + L_t m_{y_t} - \hat x_t\|^2+\|h(x_t)\|^2  }{ \|x_t  \|^2}\nn\\
&=\lambda_{max, \Theta_t}\lambda_{max, P_t}\lim_{\|x_t\|\rightarrow \infty}\frac{\|x_t \|^2+\|h(x_t)\|^2  }{ \|x_t  \|^2}\nn\\
&=\lambda_{max, \Theta_t}\lambda_{max, P_t}\left(1+\lim_{\|x_t\|\rightarrow \infty}\frac{\|h(x_t)\|^2  }{ \|x_t  \|^2}\right)\nn\\
& =\lambda_{max, \Theta_t}\lambda_{max, P_t}\left(1+\delta_h\right).\nn
 }
Since $\Theta_t\rightarrow 0$ as $c_t\rightarrow 0$, it follows that $\lambda_{max,\Theta_t}\rightarrow 0$ as $c_t\rightarrow 0$. Since $\delta_h$ is finite, we can make  
$$ \lim_{\|x_t\|\rightarrow \infty} \frac{|\mu(x_t)^\top  \Theta_t \mu(x_t)| }{ \|x_t - \hat x_t \|^2_{P^{-1}_t}}$$
arbitrarily small taking $c_t$ sufficiently small. Thus, the term $-\frac{1}{2} \|x_t - \hat x_t \|^2_{P^{-1}_t}$ dominates $\frac{1}{2}\mu(x_t)^\top  \Theta_t \mu(x_t) $ in the integral in \eqref{inttt}. Accordingly, \eqref{inttt} takes  finite value for $c_t$ sufficiently small. \qed
 \end{pf}

We conclude that $( \hat x_{t|t},\tilde{\psi}_{t}^{0})$ represents an approximate  solution to Problem (\ref{minimaxrsu2}) in the sense that:
\begin{equation}\begin{aligned}
                    J_t(\tilde \psi_t,~ \hat x_{t|t}) \leq J_t(\tilde \psi^{0}_t, \hat x_{t|t}) \approx J_t(\tilde \psi^{\star}_t, g^{\star}_{t})  \leq J_t(\tilde \psi^{\star}_t, g_t),
                    \end{aligned}\nn
                  \end{equation}
for any $\tilde \psi_t \in \mathcal{B}_t$ and $g_t \in {\mathcal{G}}_t$. Note that, the first of the above inequalities follows from Proposition  \ref{l_lfm2}. Thus, the  approximate least favorable model over the finite time interval $\tN$ takes the form:
 \begin{equation} \label{tilde_p0ZN} \begin{aligned}  \tilde p^0(Z_N) \propto  p_0(x_0) \prod_{t=0}^N   p_t(x_{t+1} | x_t) \tilde \psi^0_t(y_t | x_t) \end{aligned} \end{equation}
where the symbol $\propto$ means that the two sides are equal up to a constant factor. In what follows, we will simply refer to (\ref{tilde_p0ZN}) as least favorable model without specifying that it is an approximation.


\section{Adversarial data generation} \label{sec_lfm}
In order to assess the performance of the robust sigma point filter in the least favorable scenario, we need to develop a simulator for generating random samples from the least favorable density (\ref{tilde_p0ZN}), i.e adversarial data.  In the linear setup, it was shown that the least favorable model  admits a state space  representation over a finite time interval, \cite[Theorem 4]{URKF}. Such result, however, cannot be exploited in this nonlinear setting.  Indeed, a fundamental aspect in the linear setup is that the least favorable density  is Gaussian, but the least favorable density in (\ref{tilde_p0ZN}) is not Gaussian in general. Thus,  it is not straightforward to draw samples from (\ref{tilde_p0ZN}). Markov chain Monte Carlo (MCMC) algorithms are effective techniques for approximately sampling from complex probability densities in high-dimensional spaces. Thus, we use the Metropolis-Hastings (MH) algorithm \cite{hill2019stationarity} in order to tackle the problem. Our target density is  $\tilde p^0(Z_N)$, defined in (\ref{tilde_p0ZN}).  Suppose it is easy to generate a random sample from a proposal probability density  $ \bar q( Z_N | Y_N^k)$ where $Z_N^k=\{Y_N^k,X_{N+1}^k\}$, $Y_N^k$ and $X_{N+1}^k$ are the subvectors of $Z_N^k$ containing the  observation and state components, respectively.  We consider the MH scheme outlined in Algorithm \ref{ALGOMH}
which provides samples of  $Z_N$, following the target $\tilde p^0$ generated by the proposal $ \bar q$.  
\begin{rem}\label{remk1}
It is worth noting that both the target and proposal densities are  formed as products of multiple probability density functions, each of which includes the density of $x_0$ and the state transition density $p_t(x_{t+1}|x_t)$.  Since we know how to draw samples from $x_0$, and $x_{t+1}$ given $x_t$ (both are Gaussian distributed), we can consider  \begin{equation} \label{t_q}  \pi(Z_N):=  \prod_{t=0}^N   \tilde \psi^0_t\left(y_t | x_t\right).\end{equation}
in place of the target density. Indeed, the latter can be used to compute the acceptance ratio: \begin{align*}
\frac{\tilde p^0( Z_N)  \bar  q( Z_N^k |   Y_N )}{ \tilde p^0(Z_N^k)   \bar q( Z_N | Y_N ^k)} =
\frac{\pi( Z_N) q( Z_N^k |  Y_N )}  {\pi ( Z_N^k)    q(Z_N | Y_N^k )}
\end{align*}
where $$ q(Z_N| Y_N^k):=\frac{\bar q(Z_N| Y_N^k)}{p_0(x_0) \prod_{t=0}^N p_t(x_{t+1}|x_t)}.$$Henceforth, with some abuse of terminology, we will refer to (\ref{t_q}) as the target density.
\end{rem}
 Next, we introduce the proposal density $\bar  q(  Z_N |  Y_N^k )$, as well as how to evaluate $\pi (Z_N)$ and $q( Z_N |  Y_N^k)$. In what follows, we assume that $f:\, \Rs^{n}\rightarrow \Rs^n$ and $h:\, \Rs^{n}\rightarrow \Rs^m$ are continuous functions.

\begin{algorithm}[t]
    \caption{MH scheme for the simulator}\label{ALGOMH}
    \begin{algorithmic}[1]
        \State  Generate  $ Z^0_N$ from the nominal  model (\ref{nomi_mod})
        \For{$k\geq0$}
        \State  Draw $Z_N$ from  $\bar q(Z_N \mid Y_N^k)$ with $Y^k_N$  observations
    \NoNumber{extracted from $ Z^k_N$}
          \State  Compute the acceptance ratio
\begin{equation}
 \alpha_k = \min\left(1, \frac{\pi( Z_N)}{\pi ( Z^k_N)}  \frac{q( Z_N^{k} |    Y_N )}{  q( Z_N |Y^k_N )} \right)\nn
\end{equation}
\hspace{0.5cm} with $ Y_N$  observations extracted from $Z_N$
\State Draw $u_k$ from $\mathcal U [0,1]$
\If{$u_k \leq \alpha_k$}
\State $ Z^{k+1}_N =  Z_N$
\State $ k=k+1$
\EndIf
        \EndFor
    \end{algorithmic}
\end{algorithm}

\subsection{Proposal density} \label{sec_51}
 We construct   the proposal density relying on the least favorable state space model in  \cite[Theorem 4]{URKF} derived for the linear case.  More precisely, given the  observations $Y_{N}^k$, we can compute $\hat{x}_{t|t}$, $\hat{x}_{t+1}$, $L_t$ and $\theta_{t|t}$, with  $t\in\tN$, using  the robust sigma point Kalman filter outlined in Algorithm \ref{RUKF2}. Then, we linearize the nominal nonlinear model (\ref{nomi_mod}) along the estimated state trajectory  by  the  robust sigma point Kalman filter:
\begin{equation}\label{linearized_model}
     \begin{array}{rl}   x_{t+1} &=   A_{t}    x_t -  A_{t} \hat {x}_{t|t} +  f(\hat{x}_{t|t}) +   B v_t\\
      y_t &= C_{t}   x_t -    C_{t} \hat {x}_{t} +  h(\hat{x}_t)+ D v_t\\
 \end{array}
 \end{equation}
where  $A_t$ and $C_t$ are the Jacobian matrices:\begin{equation} \label{AC}A_{t} := \left. \frac{\partial f(x)}{\partial x}\right |_{x=\hat {  x}_{t|t}},\;  \;  C_{t}: = \left. \frac{\partial  h(x)}{\partial x}\right |_{x=\hat { x}_{t}}. \end{equation}
The least favorable  model over the time interval $\tN$ corresponding to the ambiguity set about the linearized model (\ref{linearized_model}) and with tolerance sequence $\{c_t,\; t\in\tN\}$  is given by using the result in \cite[Theorem 4]{URKF}:
\begin{equation}\label{LF_x2}\begin{aligned}
\eta_{t+1} &=\bar {A}_t \eta_{t} + \bar a_t + \bar B_t   \epsilon_t \\
y_{t} &=\bar{C}_t \eta_{t}+ b_t + \bar D_t   \epsilon_t
\end{aligned}\end{equation}
where 
\begin{equation} \label{barABCD}\begin{aligned}
& \eta_t = \left[\begin{array}{c}
                  x_t \\
                  e_{t-1} \\
                   Bv_{t-1}
                \end{array} \right],~ \bar a_t= \left[\begin{array}{c}
                   a_t \\
                 0 \\
                 0
                \end{array}\right],  ~\epsilon_t =   \left[\begin{array}{c}
                  B v_t \\
                   \upsilon_{t}
                \end{array} \right], \\           
&a_t = -   A_{t}  \hat{{x}}_{t|t} +  f(\hat{{x}}_{t|t}), \\
&\bar{A}_{t}:=\left[\begin{array}{ccc}
  A_t & 0 & 0\\
0 &  \Delta_t  - L_t F_t A_{t-1} &  \Lambda_t  \\
0 &0 &0\\
\end{array}\right],\\
&\bar {B}_{t}:=\left[\begin{array}{cc}
 I_n  & 0 \\
 0 & -L_t  \Upsilon_t\\
 I_n & 0\\
\end{array}\right],~~ b_t  = - C_{t} \hat{{x}}_{t} +  h(\hat{{x}}_t),\\
&\bar {C}_{t}:=\left[\begin{array}{llll}
C_t & F_t A_{t-1} & F_{t}
\end{array}\right], \quad \bar {D}_{t}:=\left[\begin{array}{cc}
 0 &  \Upsilon_t
\end{array}\right],\\
&\Delta_t :  = A_{t-1} - L_{t}C_t A_{t-1}, ~ \Lambda_t := I_n - L_{t} F_{t} - L_t C_t;
\end{aligned}\end{equation}
 $\varepsilon_{t} \in \Rs^{m}$ is normalized WGN; 
\begin{equation*}
\begin{aligned}e_{t} &:= \left( \Delta_t - L_t F_t  A_{t-1} \right) e_{t-1} + \Lambda_t B v_{t-1} - L_t  \Upsilon_t \varepsilon_{t}; \\
O_t&:=\left((DD^\top)^{-1} -L_t^\top(\Omega_{t+1}^{-1}+\theta_{t|t}I_n) L_t  \right)^{-1},\\
 F_{t}&  := - O_t L^\top_t (\Omega_{t+1}^{-1}+\theta_{t|t}I_n)  (I_{ n} - L_{t}  C_{t});
\end{aligned} \end{equation*}
$\Upsilon_t$ is the Cholesky factor of ${O}_{t}$, i.e. $\Upsilon_t \Upsilon^\top_t = O_t$;  $L_t$ is the Kalman filtering gain obtained by Algorithm \ref{RUKF2}.
Moreover, $\Omega^{-1}_{t}$
is computed by the following backward recursion:
\begin{equation*} \label{omgea}\begin{aligned}\Omega^{-1}_{t} = A^\top_{t-1}  F^\top_{t} O_t F_{t} A_{t-1} + \Delta_t^\top (\theta_{t|t}  I_n + \Omega^{-1}_{t+1})   \Delta_t  \end{aligned}\end{equation*}
with $\Omega^{-1}_{N+1} = 0$.  Notice that we can rewrite (\ref{LF_x2}) as
\begin{equation}
\label{linearized_lfm2}
\begin{aligned}
   x_{t+1} &=  {A_{t}  x_t }+  a_t +  B v_t \\
    y_t &= {C_{t}  x_t } +  b_t + F_t A_{t-1} e_{t-1} +F_t B v_{t-1} +\Upsilon_t\varepsilon_t 
\end{aligned}
\end{equation}
 where $e_{t-1}$ is independent from $v_{t-1}$ and it is Gaussian distributed with zero mean and covariance matrix  $\Pi_{e_t}$. The latter is the $n \times n$ submatrix of $\Pi_t$ obtained selecting rows and columns from $n+1$ to $2n$;
$\Pi_t $ is the solution to the Lyapunov equation (see \cite[Section 3]{URKF}):
 \al{\Pi_{t+1}=\Gamma_t \Pi_t \Gamma_t^\top +X_t\Xi X_t^\top \nn}
 where
{\small \begin{align*} \label{a_lfm1}
& \Gamma_t:=\left[\begin{array}{ccc}
\Delta_t &  -L_{t} F_{t} A_{t-1} &\Lambda_t \\
0 & \Delta_t - L_t F_{t} A_{t-1} & \Lambda_t \\
0 & 0& 0
\end{array}\right], \;
X_t=\left[\begin{array}{cc}
0  & -L_t  \Upsilon_t\\
 0   & -L_t  \Upsilon_t\\
  I_n & 0
\end{array}\right],
\end{align*}}and
$\Xi$ is the block diagonal matrix with main blocks $\{ B B^\top,I_m\}$. 
Therefore, the measurement model in (\ref{linearized_lfm2}) is characterized by the transition  density of $ y_t$ given $x_t$ and $Y_N^k$ which is Gaussian with mean $ C_{t} x_t + b_t$ and covariance matrix $  Q^L_t $
where  \begin{align*}  Q^L_{t} &:= F_{t} A_{t-1} \Pi_{e_{t-1}} A^\top_{t-1} F^{\top}_{t} + F_{t} B B^\top F^\top_{t} + \Upsilon_t \Upsilon^\top_t.\end{align*}
However, our aim is to  develop a proposal density that effectively captures the essential characteristics  of the target density.  Notice that, $ A_t x_t +a_t $ and  $ C_{t} x_t + b_t$  represent the first order Taylor expansion of $f(x_t)$ and $h(x_t)$  around $\hat{ x}_{t|t}$ and $\hat{ x}_{t}$, respectively. Then, a refined version of (\ref{linearized_lfm2}) is the one in which $ A_t x_t +a_t $ and  $ C_{t} x_t + b_t$  are replaced by $f(x_t)$ and $h(x_t)$ , respectively.
The corresponding ``proposal'' state space model is
\begin{equation}\label{linearized_lfm3}
     \begin{array}{rl}    x_{t+1} &=  f(x_t) +    B v_t  \\
    y_t &= h(x_t) + F_t A_{t-1} e_{t-1} +F_t B v_{t-1} +\Upsilon_t\upsilon_t \\
\end{array}
\end{equation}
and the corresponding proposal density takes the form
\begin{equation} \label{defPROP} \bar q( Z_N|Y_N^k)= p_0( x_0)\prod_{t=0}^N p_t( x_{t+1}|x_t)\psi^L_t\left({  y_t | x_t ,Y_{t-1}^k}\right) \end{equation}
 where   \begin{equation}\label{psiL}{\psi^L_t( y_t| x_t,  Y_{t-1}^k)}  = \mathcal N (h(x_t), Q^L_{t}).\end{equation}
To sum up, the state space model (\ref{linearized_lfm3}) is constructed as follows: first,  given  the  observations $Y_N^k$, we perform a forward sweep to compute  $\hat x_{t|t}$, $\hat x_{t+1}$, $L_t$ and $\theta_{t|t}$ by means of the proposed robust sigma point Kalman filter, then we perform a backward sweep to compute $\Omega^{-1}_t$ and finally we perform a forward sweep to compute the matrices characterizing the state space model. In this respect, it is easy to draw a proposal  sample $Z_N$ from the state space model (\ref{linearized_lfm3}). Finally,   in view of Remark \ref{remk1},  we only need to evaluate
\begin{equation} \label{bar_p2}  q({Z_N |  Y_N^k})=  \prod_{t=0}^N  \psi^L_t\left({  y_t |x_t ,Y_t^k}\right).\end{equation}
Accordingly, in view of (\ref{psiL}) and (\ref{bar_p2}), we obtain
\begin{equation}
\begin{aligned}
 q( Z_N | Y_N^k ) =  \prod_{t=0}^{N}\frac{1}{ \bar M_t} \exp\left( - \frac{1}{2}\|y_t - h(x_t) \|^2_{( Q^L_t)^{-1}} \right)
 \end{aligned}\nn
\end{equation}
where  $\bar M_t = \sqrt{(2 \pi)^{m}\left| Q^L_t\right|}.$

\subsection{Target density} \label{sec_52}
 In regard to the evaluation of  target density $\pi(Z_N) $,  by  (\ref{psi0_t}) and (\ref{t_q}), we have 
  \begin{align}
\pi(Z_N)  =  \prod_{t=0}^N \frac{1}{M_{t}} \exp \left(\frac{ \theta_{t|t}}{2} \left\|  x_{t}- \hat x_{t|t}\right\|^{2}\right) \psi_{t}(y_t|x_t),\nn\end{align}
 so the challenging aspect regards  the computation of $M_t$ defined in (\ref{M_T3}).   It is worth noting that it is not possible to find a closed form expression
 for the latter  because of the presence of $\psi_t(y_t|x_t)$, i.e. its mean is a nonlinear function of $x_t$. Then,  an approximation of $M_t$ can be obtained through Monte Carlo integration:
 \begin{equation*}\label{t_M2}  \hat M_{t,r}:=   \frac{1}{r} \sum^{r}_{j=1} N_t (x_t^j),  \end{equation*}
where $ N_t(x_t)$ has been defined in \eqref{Nxt2}; $ x^{1}_t \ldots  x^r_t$  are independent samples drawn from  $\tilde {p}_t^a(x_t| Y_{t-1}) = \mathcal N(\hat x_t, P_t)$.  Note that, $N_t(x_t)$ is a random variable obtained from a nonlinear transformation of $x_t$ given $Y_{t-1}$. The next theorem shows under which conditions $N_t(x_t)$ has finite variance.
 \begin{prop} Let $\sigma^2_t$ denote the variance of $N_{t}(x_t)$.
  Assume that $h(x)$ is such that  $\delta_h$ defined in \eqref{h_cond} is finite.
  Taking the tolerance $c_{t}$ sufficiently small, then $\sigma_t^2$ is finite and 
  \al{\label{clt} \sqrt{r}(\hat M_{t,r}-M_t)\overset{p}{\longrightarrow} \Nc(0,\sigma^2_t).}
\label{propos2minor}
\end{prop}
\begin{pf} 
 First, note that  
\[
\sigma^2_t = \mathbb{E}[N^2_t(x_t)] - (\mathbb{E}[N_t(x_t)])^2 \geq 0.
\]
Thus, it suffices to prove that $ \mathbb{E}[N^2_t(x_t)]$ is finite.  We have
\al{\label{defEN}\begin{aligned} \mathbb{E}[N^2_t&(x_t)]  = \int N^2_t(x_t) \tilde p_t^a(x_t|Y_{t-1})dx_t  \\
&=  \frac{1}{|S_t||R|}\int  \frac{1}{\sqrt{(2 \pi)^n |P_t|}} \exp(\mu(x_t)^\top \Theta_t \mu(x_t) \\
 & \hspace {3cm}- \frac{1}{2} \|x_t - \hat x_t \|^2_{P^{-1}_t})dx_t .
\end{aligned} }
Since $\delta_h$ is finite and $c_t$ is taken sufficiently small (using the same arguments in the proof of Theorem \ref{propos2}) we have that  the term $-\frac{1}{2} \|x_t - \hat x_t \|^2_{P^{-1}_t}$ dominates $\mu(x_t)^\top  \Theta_t \mu(x_t) $ in the integral. Accordingly, \eqref{defEN} takes finite value for $c_t$ sufficiently small.  Since $N_t(x_t)$ has finite variance and the samples $N_t(x_t^1)\ldots N_t(x_t^r)$ are independent, by the central limit theorem we obtain \eqref{clt}. \qed
 \end{pf}
 By Proposition \ref{propos2minor}, we have that the following condition approximately holds for $r$ large 
$$\hat M_{t,r} \sim \mathcal N\left(M_{t},~ \frac{\hat \sigma^2_{t,r}}{ r} \right)$$
and 
 \begin{equation}
       \hat \sigma^2_{t,r} := \frac{1}{r-1} \sum_{j=1}^{r}(N_t(x^j_t) - \hat M_{t,r})^2\nn
  \end{equation}
is the sample variance estimator of $\sigma^2_t$.
Thus,  it remains to design the number of samples $r$ in such a way that $\hat M_{t,r}$ approximates $M_t$ with a certain accuracy.  
With a $95\%$ confidence level, it holds that
$$ \frac{\sqrt{r}|\hat M_{t,r} - M_t|}{\hat \sigma_{t,r}} \leq 1.96,$$
and thus the following inequality approximately holds $$\frac{|\hat M_{t,r} - M_t| }{M_t}\leq \frac{1.96\,\hat  \sigma_{t,r}}{\sqrt{r}\hat M_{t,r}}:=\tau_r.$$ Therefore, given the desired relative accuracy $\tau^\star$ on the computation of $M_t$ and an initial value for $r$, we check whether $\tau_r\leq \tau^\star.$ If not, we increase $r$ until $\tau_r\leq \tau^\star.$
Notice that, if we need to increase $r$, then $\hat  \sigma_{t,r}$ and $\hat M_{t,r}$ can be updated recursively using well known formulas.

 In the  special case $h(x)$ is a linear function, i.e. $h(x) = C x$, the normalization constant $M_t$ can be  computed in closed form. 

 \begin{prop} \label{prophlin}Let $h(x)=Cx$. Taking the tolerance $c_t$ sufficiently small, then
\begin{align} \label{MT5}
M_t = & \frac{\exp(\frac{1}{2} (   \|\hat x_t - L_t m_{y_t} \|^2_{\Gamma_t}-\|\hat x_t \|^2_{P_t^{-1}} + \|\zeta_t\|^2_{U_t^{-1}}  )}{\sqrt{|S_t||R||P_t||U_t|^{-1}}}\nn
\end{align}
where 
\begin{align}
\Gamma_t=& \theta_{t|t}I+\theta_{t|t}^2 L_t S_t^{-1} L_t^\top\nn\\
\zeta_t =& \theta_{t|t} U_t  ((C^\top R^{-1} - \theta_{t|t} L_t) S_t^{-1} L_t^\top - I)    (\hat x_t - L_t m_{y_t} ) \nn\\
&  \hspace{.2cm} + P_t^{-1} \hat x_t \nn \\
U_t = & (C^\top R^{-1} C + P_t^{-1} - \theta_{t|t} I \nn \\
&\hspace{.2cm}  -(C^\top R^{-1} - \theta_{t|t} L_t) S_t^{-1} (R^{-1} C - \theta_{t|t} L_t^{\top}))^{-1} 
\end{align} 
and $R$, $S_t$ have been already defined in \eqref{slSR}.
 \end{prop} 
 
 \begin{pf}
 In view of (\ref{M_T3}),  (\ref{Nxt2}), (\ref{hat_M_2}) and $\tilde p_t^a(x_t | Y_{t-1}) =\mathcal N(\hat x_t, P_t)$, we have
 \begin{align}
 M_t & = \int N_t(x_t) \tilde p^a_t(x_t | Y_{t-1}) d x_t \nn\\
 &\hspace{-0.3cm} =  \frac{\int \exp(\frac{1}{2} (s_t(x_t)^\top S_t^{-1} s_t(x_t) + l_t(x_t) - \| x_t - \hat x_t \|^2_{P_t^{-1}} )) d x_t}{\sqrt{(2 \pi)^n |S_t||R||P_t| }} \nn
\end{align} 
where 
\al{ l_t(x_t) &= \theta_{t|t} \| x_t - \hat x_t + L_t m_{y_t} \|^2 - \|  Cx_t\|^2_{R^{-1}},\nn\\
s_t (x_t)&= \theta_{t|t} L_t^\top (\hat x_t - x_t - L_t m_{y_t}) + R^{-1} Cx_t.\nn } 
By \eqref{defSi}, we have 
\al{ C^\top & R^{-1} C + P_t^{-1} - \theta_{t|t} I \nn \\
&\hspace{.5cm}  -(C^\top R^{-1} - \theta_{t|t} L_t) S_t^{-1} (R^{-1} C - \theta_{t|t} L_t^{\top})\nn\\
&= C^\top R^{-1} C + P_t^{-1} - \theta_{t|t} I-(C^\top R^{-1} - \theta_{t|t} L_t) \nn \\
&\hspace{.5cm}  \times (R+\theta_{t|t} RL_t^\top L_t R+\mathcal{O}(\theta_{t|t}))(R^{-1} C - \theta_{t|t} L_t^{\top})\nn\\
&=C^\top R^{-1} C + P_t^{-1}-C^\top R^{-1} C +\mathcal{O}(1) \nn\\ &=P_t^{-1}+\mathcal{O}(1)\nn}
where $\mathcal{O}(1)$ are the infinitesimal terms of order higher than or equal to $\theta_{t|t}$. Thus, $U_t$ is positive definite (and thus also invertible) for $c_t$ taken sufficiently small (recall that $\theta_{t|t}\rightarrow 0$ as $c_t\rightarrow 0$, see the proof of Theorem \ref{propos2}).
After some linear algebra manipulations, we have
\al{&s_t(x_t)^\top  S_t^{-1} s_t(x_t) + l_t(x_t) - \| x_t - \hat x_t \|^2_{P_t^{-1}}\nn\\
&\hspace{0.5cm}= -\| x_t-\zeta_t\|^2_{U_t^{-1}}+\|\hat x_t-L_t m_{y_t}\|^2_{\Gamma_t}- \|\hat x_t\|^2_{P_t}+\|\zeta_t\|^2_{U_t^{-1}}.
\nn }
Thus, {\small \begin{equation}
\begin{aligned}
 M_t &=   \frac{\exp(\frac{1}{2} (   \|\hat x_t - L_t m_{y_t} \|^2_{\Gamma_t}-\|\hat x_t \|^2_{P_t^{-1}} + \|\zeta_t\|^2_{U_t^{-1}}  )  }{\sqrt{(2 \pi)^n |S_t||R||P_t| }}\nn\\
 & \hspace{0.5cm}\times  \int \exp(-\frac{1}{2}  \| x_t - \zeta_t \|^2_{U^{-1}_t} )  d x_t\nn\\
 &\hspace{-0.5cm}= \frac{\exp(\frac{1}{2} (   \|\hat x_t - L_t m_{y_t} \|^2_{\Gamma_t}-\|\hat x_t \|^2_{P_t^{-1}} + \|\zeta_t\|^2_{U_t^{-1}}  )  }{\sqrt{(2 \pi)^n |S_t||R||P_t| }} \sqrt{(2\pi)^n |U_t |}\nn\\
\end{aligned}
\end{equation}}which proves the claim.\qed
 \end{pf}

\begin{rem} Although the measurement model in Proposition \ref{prophlin} is linear, this does not lead to a linear filtering problem such as the one in \cite{URKF}, since the process model remains nonlinear.
 \end{rem}
Finally, the next results establish the convergence properties of the proposed algorithm.
\begin{thm} Assume that $\delta_h$ defined in \eqref{h_cond} is finite. Taking the tolerance $c_t$ sufficiently small, the Markov chain $Z_N^k$ generated by Algorithm \ref{ALGOMH} converges in total variation distance to the target density $\tilde p^0$ defined in \eqref{tilde_p0ZN}, i.e. for every $Z^0_N$ and $\varepsilon>0$ there exists a positive  integer $k^\star$ such that for every set $\mathcal A\subseteq \Rs^{(N+2)n+(N+1)m} $
  \begin{equation}
    \left|\mathbf{P}[Z^k_N \in \mathcal A | Z^0_N ] - \int_{\mathcal A} \tilde p^0(Z_N) d Z_N\right|<\varepsilon \hbox{ for } k\geq k^\star\nn
  \end{equation}
  where $\mathbf{P}[Z^k_N \in \mathcal A | Z^0_N ]$ is the probability that $Z^k_N \in \mathcal A$ if the initial condition is equal to $Z^0_N$.  
  \end{thm} 
  \begin{pf} The target \eqref{tilde_p0ZN} and the proposal \eqref{defPROP} are strictly positive functions, since they are given by products of exponential terms and positive normalizing constants. Hence:
  \begin{itemize}
  \item For any hypercube $\mathcal T\subseteq  \Rs^{(N+2)n+(N+1)m}$
  \al{&\inf_{Z_N, \tilde Z_N\in \mathcal T}\left[\bar q( \tilde Z_N |    Y_N )  \min\left(1, \frac{\tilde p^0( Z_N)}{\tilde p^0 ( \tilde Z_N)}  \frac{\bar q( \tilde Z_N |    Y_N )}{  \bar q( Z_N |\tilde Y_N )} \right) \right]>0\nn\\
  & \int _{\mathcal T}\tilde p^0(Z_N) d Z_N>0\nn}
  where $Z_N=\{ Y_N, X_{N+1}\}$ and $\tilde Z_N=\{\tilde Y_N,\tilde X_{N+1}\}.$
   \item Any point in $\Rs^{(N+2)n+(N+1)m}$ can be connected to any subset $\mathcal T\subseteq \Rs^{(N+2)n+(N+1)m}$, such that $\int_{\mathcal T}\tilde p^0(Z_N) dZ_N>0$, using a finite number of balls that have radii less than or equal to $\eta/2$ and overlap in sets over which the integral of  $\tilde p^0$ is positive; moreover,
  $$  \bar q(Z_N| \tilde Y_N)>0 \hbox{ for any } \| Z_N-\tilde Z_N \|<\eta \hbox{ with } \eta>0.$$ 
 Thus, the transition kernel is irreducible with respect to the target.
\end{itemize}
   Accordingly, by  \cite[Conditions C1$^\prime$ and C2$^\prime$]{hill2019stationarity}, the convergence in total variation distance holds.\qed 
  \end{pf}

\section{Numerical Experiments}\label{sec_8}
\subsection{ Worst-case performance} 
  \begin{figure}
  \centering
  \includegraphics[width=0.5\textwidth]{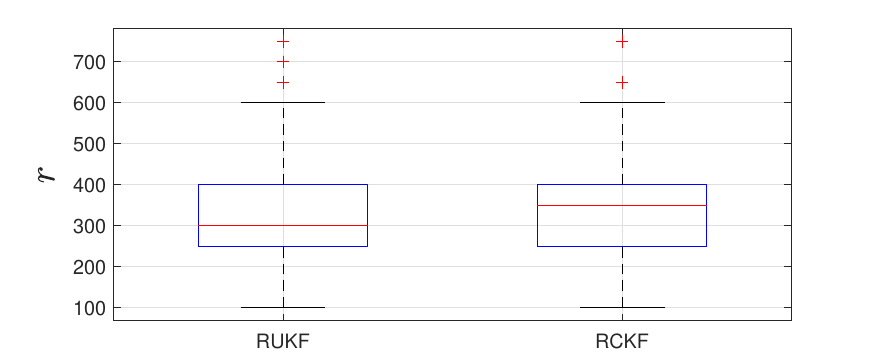}
  \caption{Boxplot of the value of $r$ needed in the MH scheme  and corresponding to the least favorable model of  RUKF and RCKF.}\label{FIG_r}
\end{figure}

 \begin{figure}
  \centering
  \includegraphics[width=0.5\textwidth]{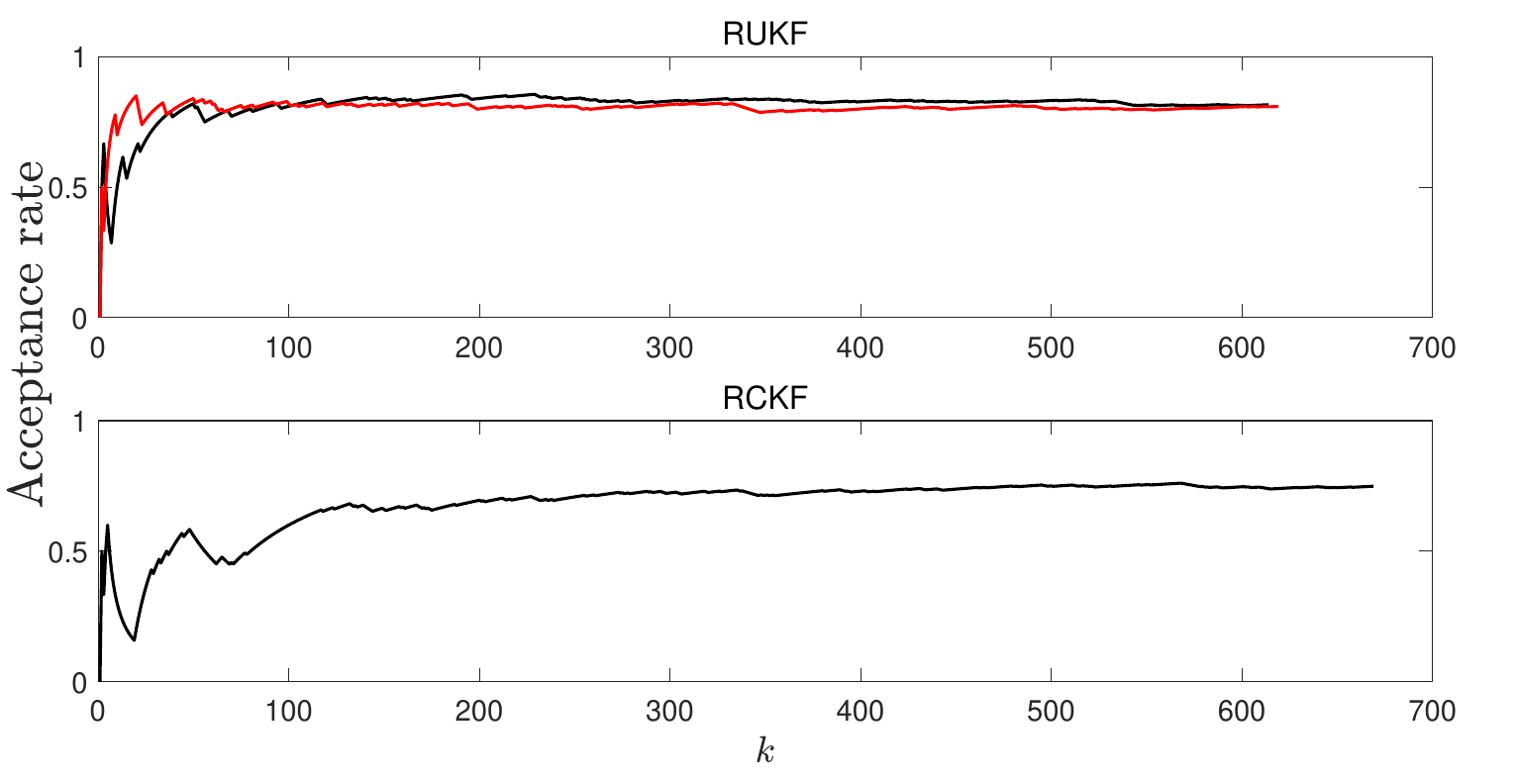}
  \caption{Acceptance rate in the MH algorithm corresponding to the least favorable model of RUKF and RCKF. Variable $k$ denotes the number of proposals (including both accepted and rejected proposals). For the RUKF case, we have depicted the trajectories corresponding to two different samples of $Z_N^0$.}\label{FIG_acc}
\end{figure}

We  assess the performance of the proposed filters in the worst-case scenario (i.e. using adversarial data).  We consider the nominal  nonlinear state space model \eqref{nomi_mod}  with 
\[
\left\{
\begin{array}{rcl}
f(x) & = &
\begin{bmatrix}
0.2x_1+x_2+0.5\cos(2x_2)-1\\
0.95x_2
\end{bmatrix}
\\[1.2em]
h(x) & = &
\begin{bmatrix}
x_1+0.1\left(\sqrt{1+x_2^2}-1\right)\\
0.1x_1-0.1x_2
\end{bmatrix}.
\end{array}
\right.
\]
 Moreover,
$x=[\,  {x}_{1}^{\top} \; {x}_{2}^\top \,]^\top$, $x_0\sim \mathcal N ( 0,  I_2)$ and
$B=[\, I_2 \; 0\,]$,  $D=[\,0  \; I_2 \,]$.
Notice that, the chosen nominal model satisfies the hypothesis on $h$  in Theorem 2.
In what follows, we consider the following nonlinear filters:
\begin{itemize} 
  \item UKF denotes the unscented Kalman filter  \cite{julier2000new} where the parameters in (\ref{W_ukf}) are set as \( a = 0.5 \), \( b = 2 \), and \( \kappa = 1 \);
   \item CKF denotes the cubature Kalman filter \cite{arasaratnam2009cubature};
  \item RUKF  denotes the robust sigma point Kalman filter of Section \ref{sec_3_true} with tolerance \( c =2 \cdot 10^{-4} \); the sigma points are obtained using the unscented transformation whose parameter setting is the same as the one of UKF;
  \item RCKF  denotes the robust sigma point Kalman filter of Section \ref{sec_3_true} with tolerance \( c = 2 \cdot 10^{-4} \); the sigma points are obtained using the spherical cubature rule.
\end{itemize}
We evaluate their performance under the least favorable models corresponding to RUKF and RCKF introduced above.
More precisely,  we have generated $M=500$ adversarial samples $Z^k_N$ of length $N=50$ from each least favorable model by means of the MH scheme proposed in Section \ref{sec_lfm}. We set the initial value of \( r \) equal to \( 100 \) and the desired relative accuracy  \( \tau^\star = 2 \cdot 10^{-3} \) for the computation of \( \hat{M}_{r,t} \). In this way, the cumulative relative error over the time interval  $[0,50]_{\mathbb Z}$ is approximately equal to 10\%. An upper bound equal to 1000 is imposed on $ r$ to control the computational time in edge cases. Fig. \ref{FIG_r} shows the boxplot of the value of $r$ required to compute $\hat M_{r,t}$  (both rejected and accepted samples have been considered) by the two MH algorithms.  We can see that the typical range of $r$ is 250$\div$400 and values greater than 600 are considered as outliers; thus, the chosen upper bound for $r$ is adequate.  Fig. \ref{FIG_acc} shows the corresponding acceptance rate for the  two different simulators.
 Regarding the least favorable model of RUKF, i.e. the first subfigure in Fig. \ref{FIG_acc},  we have considered the  trajectories generated by two different  samples of $Z^0_N$. Notably, the algorithm exhibits a similar burn-in period and acceptance rate for the two trajectories.  Regarding the least favorable model of RCKF, as in the second  subfigure in Fig. \ref{FIG_acc}, we have depicted  only one trajectory, as we observed a similar behavior.  Overall, the acceptance rate in both cases  is approximately equal to  \(80\%\),  indicating that the proposal density has been properly designed and provides a good approximation to the target density. 
  \begin{figure}[t]
  \centering
  \includegraphics[width=0.5\textwidth]{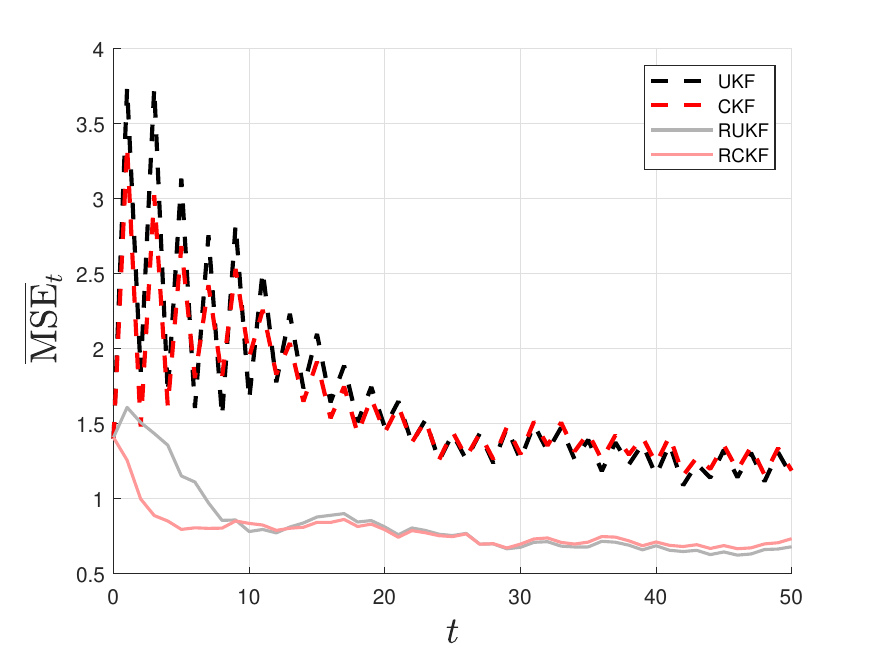}
  \caption{Mean squared error of the filters using adversarial data corresponding to the least favorable model of   RUKF.}\label{FIG_RUKF}
\end{figure}
\begin{figure}[t]
  \centering
  \includegraphics[width=0.5\textwidth]{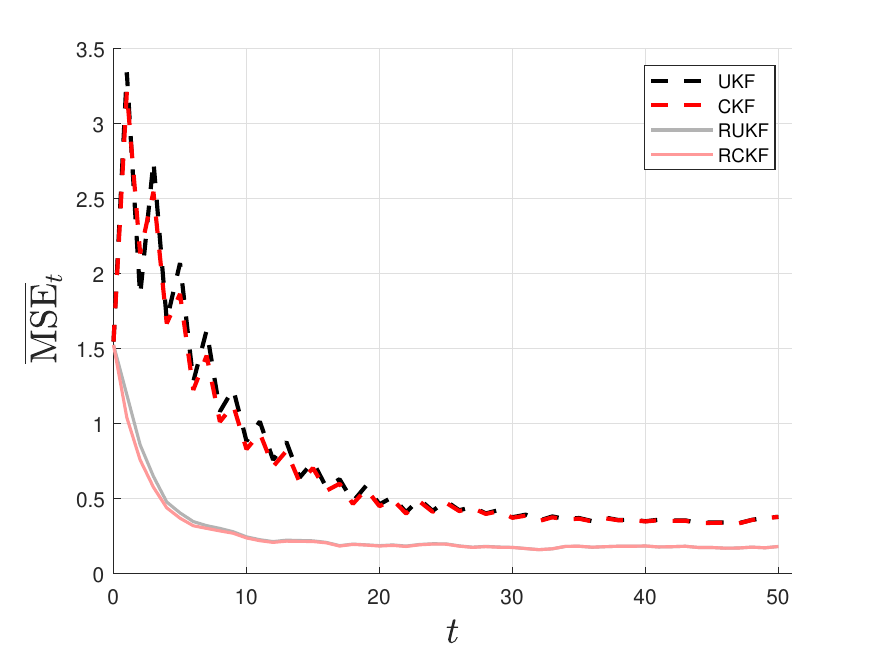}
  \caption{Mean squared error of the filters using adversarial data corresponding to the least favorable model of   RCKF.}\label{FIG_RCKF-eps-converted-to.pdf}
\end{figure} Next,  we  assess the filters using the aforementioned dataset. More precisely, for each sample $Z_N^k$ we extract the state trajectory  $X_N^k$ and the corresponding measurement trajectory  $Y_N^k$. Then, we evaluate the performance of the  filters computing the mean squared error at time $t$: 
\begin{equation} \label{def_perf}
  \overline{\text{MSE}}_t =  \frac{1}{M} \sum^M_{k=1} \| {{x}}_t^k-\hat{x}_{t|t}^k \|^2
\end{equation}
where $x_t^k$ is the state at time $t$ extracted from $X_N^k$, while  $\hat{x}^k_{t|t}$ is the filtered state estimate obtained by the observations $Y^{k}_N$. Fig. \ref{FIG_RUKF}  and Fig. \ref{FIG_RCKF} show the mean squared error  for the four filters using the adversarial data generated by the least favorable model of RUKF and RCKF, respectively.
The results highlight that, in each case, the robust filters outperform the standard filters. Moreover, the robust filters converge to a stable regime more rapidly, while also avoiding the pronounced fluctuations observed during the initial stage.

Finally,  we consider a nominal model that is identical to the previous one,  except that the state transition function $f(x)$, which is now chosen as
\begin{equation}
f(x)=
\begin{bmatrix}
0.2x_1+x_2+0.5\cos(x_2)-1\\
1.2x_2
\end{bmatrix}.
\end{equation}
Under adversarial data, the state trajectory sometimes  has a drift. In this situation, we observed that the standard filters completely lose the ability to effectively track the true state, whereas robust filters can still maintain satisfactory tracking performance.
More precisely, we evaluate the filter’s performance using the average mean squared error:
$$
\overline{\mathrm{MSE}}
=
\frac{1}{N}\sum_{t=1}^{N}\overline{\mathrm{MSE}}_t.
$$
The results are reported in Table~\ref{tab:unstable_case}.
This finding is particularly meaningful in practical applications, since the presence of a drift may arise in scenarios such as tracking and navigation, where model uncertainty is often unavoidable.

\begin{table}[tbp]
\centering
\caption{$\overline{\mathrm{MSE}}$ using adversarial data corresponding to the least favorable model of RUKF (LF-RUKF) and RCKF (LF-RCKF).}
\label{tab:unstable_case}
\begin{tabular}{lcccc}
\hline
 & UKF & CKF & RUKF & RCKF \\
\hline
LF-RCKF & 8651900 & 8651030 & 0.8026 & 0.8017 \\
LF- RUKF & 6972663 & 6973396 & 0.8284 & 0.8283 \\
\hline
\end{tabular}
\end{table}

\subsection{Mass-spring system}
\begin{figure}[htbp]
\centering
\begin{tikzpicture}[scale=0.8]

  \fill[pattern=north east lines] (-1.2,-1) rectangle (5,-0.7);

  \draw[thick] (-1,1.2) -- (-1,-1);
  \fill[pattern=north east lines] (-1.2,-1) rectangle (-1,1.2);

  \draw[thick] (-1,0) -- (-0.8,0);
  \draw[thick,decorate,decoration={zigzag, segment length=4, amplitude=3}] (-0.8,0) -- (0.5,0);
  \draw[thick] (0.5,0) -- (1,0);

  \draw[thick,fill=white] (1,0.6) rectangle (2.5,-0.7);
  \node at (1.75,0) {$m$};

  \draw[->,thick] (2.25,0.8) -- (1.55,0.8) node[above] {$F_f$};

  \draw[->,thick] (0,0.3) -- (-0.5,0.3) node[above] {$F_{s}$};

  \draw[->,thick] (2.5,0) -- (3.3,0) node[right] {$F_e$};

  \draw[->,thick] (1.75,-0.4) -- (3,-0.4) node[right] {$p$};

\end{tikzpicture}
\caption{Mass-spring system.}
\label{fig:block-spring}
\end{figure}

In what follows, we assess the effectiveness of the proposed filters in a more challenging and realistic setting: the true model is not guaranteed to belong to the ambiguity set. We consider the mass-spring system depicted in Fig.~\ref{fig:block-spring}. The reference position of the mass is chosen such that the restoring force is equal to zero. 

{\em  Actual state space model.} Let $x=[p~s]^\top$ be the state vector, where $p~[\texttt{m}]$ is the displacement and $s~[\texttt{m/s}]$ is the velocity of the mass. The true nonlinear process model is described by the dynamics of a mass-spring system subject to friction, spring restoring force, and an external   force. Such  model can be written as  (\ref{nomi_mod}) with
\begin{equation}\label{state_space_mass_spr_m}
\begin{aligned}
  f(x)  = x +\frac{T_s}{m}  \left[\begin{array}{c}
m s \\
-F_{f} - F_{s}
\end{array}\right], ~
 B  = T_s  \left[\begin{array}{ccc}
 \epsilon & 0 & 0 \\
 0 &   \sqrt q  / m & 0
\end{array}\right]
\end{aligned}
\end{equation}
where $m=1~\texttt{[kg]}$ is the  wight of the  mass,  $T_s=0.1~\texttt{[s]}$ is the sampling time, $\epsilon = 0$ so that the relationship between displacement and velocity is  exact. Specifically, $F_f~[\texttt{N}]$ denotes the resistive force due to friction, $F_s~[\texttt{N}]$ represents the restoring force of the spring; $F_e~[\texttt{N}]$ is the external force, which is modeled as WGN with zero mean and variance $q=0.25$,   thus the second component of $Bv_t$ represents $T_sF_e/m$. 
The resistive force consists of static, Coulomb, and viscous friction components~\cite{khalil2002nonlinear}, i.e.,
\begin{equation}\label{F_f}
F_f = \alpha s + \eta(x),
\end{equation}
where $\alpha=0.5~\texttt{[Ns/m]}$ is the viscous friction coefficient, and $\eta(x)$ is defined by the following piecewise function:
\begin{equation*}\label{friction_f}
\eta(x)=
\left\{
\begin{array}{ll}
 \mu_k m g \operatorname{sign}(s),  
 & \text{for } |s|>0, \\[1mm]
 -k p,   
 & \text{for } s=0 \text{ and } |p|\leq \mu_s m g/k, \\[1mm]
 -\mu_s m g \operatorname{sign}(p),  
 & \text{for } s=0 \text{ and } |p|>\mu_s m g/k,
\end{array}
\right.
\end{equation*}
where $g=9.81~[\texttt{m/s}^2]$ is the gravitational acceleration. The parameters $\mu_k$ and $\mu_s$ denote the kinetic and static friction coefficients, respectively, which are given by
\[
\mu_k \sim \mathcal{U}(0.3,0.8), ~~~
\mu_s \sim \mathcal{U}(0.3,0.8).
\]
The restoring force is modeled by a hardening spring~\cite{khalil2002nonlinear}, where a small displacement increment beyond a certain threshold leads to a large increase in force, i.e.,
\begin{equation}\label{spring_f}
F_s = kp + k a^2 p^3,
\end{equation}
where $k=10~\texttt{[N/m]}$ is the spring constant, and $a \sim \mathcal{U}(0.01,0.05)$ is the nonlinear stiffness parameter characterizing the hardening effect.  We assume the displacement $p$ is measured by a sensor. More precisely,  four different measurement non-idealities are considered:\\
$\bullet$  sensor drift 
\begin{equation*}
y_t = p_t + \tilde{\varepsilon}_t, \qquad
\tilde{\varepsilon}_t \sim \mathcal{N}(0.5,1);
\end{equation*}
$\bullet$  intermittent observations 
\begin{equation*}
y_t = \gamma_k p_t + \tilde{\varepsilon}_t, \qquad
\tilde{\varepsilon}_t \sim
 \mathcal{N}(0,1)
\end{equation*}
where $\gamma_k$ is a binary random variable taking values in  $ \{0,~1 \}$ with success probability $\mathbb{P}[\gamma_{i,k}=1] = 0.8$;\\
$\bullet$ outlier-contaminated noise 
\begin{equation*}
y_t = p_t + \tilde{\varepsilon}_t, \qquad
\tilde{\varepsilon}_t \sim
\begin{cases}
\mathcal{N}(0, 1), & \text{with probability } 0.9, \\
\mathcal{N}(0, 5), & \text{with probability } 0.1;
\end{cases}
\end{equation*}
$\bullet$ nonlinearity dead zone 
\begin{equation*}
y_t =
\begin{cases}
p_t + \tilde{\varepsilon}_t, 
& \text{if } |p_t + \tilde{\varepsilon}_t| \ge 0.1, \\
0, 
& \text{if } |p_t + \tilde{\varepsilon}_t| < 0.1,
\end{cases}
\end{equation*}
where $\tilde{\varepsilon}_t \sim \mathcal{N}(0,1)$.

{\em Monte Carlo experiments.} We consider Monte Carlo experiments with $M=1000$ independent trials. 
In each trial, the actual process model is specified by independently sampling the uncertain parameters \(a\), \(\mu_k\), and \(\mu_s\) from the aforementioned uniform distributions. In addition, the corresponding measurement  data are generated according to one of the four  actual measurement models described above. Each trajectory is simulated over \(N=100\) time steps, with the initial state given by
$x_0 \sim \mathcal{N}\left([3\;0]^\top,\,0.1I\right)$.
Since the sampling time is \(T_s=0.1\,\mathrm{s}\), this corresponds to a total simulation duration of \(10\) seconds.  For each generated observation, we estimate the state using as nominal model 
(\ref{state_space_mass_spr_m}) with nominal parameters
\[
a=0.03, \qquad \mu_k=0.6, \qquad \mu_s=0.5, \qquad \epsilon=10^{-8},
\]
where  a small  positive  value of  $\epsilon$ is introduced  to  ensure the invertibility of the matrix $BB^\top$, and nominal measurement model 
\[
y_t = p_t + D v_t,
\]
where $D=[\, 0 \; 0 \;  1\,].$ To this end,
\begin{table*}[t]
\caption{%
$\overline{\text{MSE}}$ comparison for various filters under different measurement non-idealities.
RUKF, RCKF, and REKF reduce to the standard UKF, CKF, and EKF, respectively, when $c=0$.
}
\label{tab:rmse_comparison}
\centering
\begin{subtable}{0.48\textwidth}
\centering
\caption{Sensor drift}
\adjustbox{max width=\linewidth}{
\begin{tabular}{c|ccc|cc|cc|cc}
\toprule
$c$ & RUKF & RCKF & REKF & $N_p$ & PF & $\tau$ & ME-PF & $\sigma$ & MC-UKF \\
\midrule
0      & 3.9484 & 4.8946 & 8.5529 & 100   & 7.7581 & 0.0001 &    9.7622     & 1  & 12.661 \\
0.0001 & 3.7535 & 4.0572 & 6.1566 & 500   & 7.1670 & 0.01 &  9.3338    & 2  & 11.577 \\
0.001  & 3.3516 & 2.7624 & 42.419 & 1000  & 6.9366 & 0.03& 9.0117  & 4  & 10.690 \\
0.01   & 2.6438 & 1.5755 & 59.967 & 2000  & 6.7567 & 0.05& 8.8956 & 6  & 10.438 \\
0.03   & \textbf{2.4385} & \textbf{1.4803} & 110.04 & 5000  & 6.5137 & 0.1 & 8.7962  & 8  & 10.235 \\
0.05   & 2.4923 & 1.4993 & \textbf{3.7720} & 10000 & 6.4208 & 0.3&  8.3627   & 10 & \textbf{10.191} \\
0.1    & 2.7096 & 1.7602 & 3.9001 & 20000 & 6.3379 & 0.5& \textbf{8.2880}    & 15 & 10.757 \\
0.3    & 4.0017 & 2.0895 & 4.5181 & 50000 &      \textbf{6.0099}   &  0.8 & 8.3549& 20 & 12.401 \\
\bottomrule
\end{tabular}}
\end{subtable}
\hfill
\begin{subtable}{0.48\textwidth}
\centering
\caption{Intermittent observations}
\adjustbox{max width=\linewidth}{
\begin{tabular}{c|ccc|cc|cc|cc}
\toprule
$c$ & RUKF & RCKF & REKF & $N_p$ & PF & $\tau$ & ME-PF & $\sigma$ & MC-UKF \\
\midrule
0      & 5.6829 & 8.132 &11.4466 & 100   & 9.2916 & 0.0001 & 10.644     & 1  & 13.63 \\
0.0001 & 5.492& 7.4874 & 9.7861 & 500   & 9.0354 &  0.01 &10.253  & 2  & 12.216 \\
0.001  &5.0716 & 6.0319 & 10.468 & 1000  & 8.9495 & 0.03 & 10.091 & 4  & 11.163 \\
0.01   & 4.1387 & 3.788 & \textbf{7.8918} & 2000  & 8.898 & 0.05 &  10.464 & 6  &  10.646 \\
0.03   & 3.5296 & 2.9734 & 470.87 & 5000  & 8.7662 & 0.1  &10.561   & 8  & 10.861 \\
0.05   & \textbf{3.3623} & 2.7083 & 17.469 & 10000 & 8.6641 & 0.3 & 10.161  & 10 & \textbf{10.778} \\
0.1    & 3.3892 & 2.4101 & 11.861 & 20000 & 8.5722 & 0.5 & 10.032   & 15 & 11.089 \\
0.3    & 4.3961 & \textbf{2.132}  &1663.4 & 50000 &    \textbf{8.3237}    & 0.8  & \textbf{9.9637} & 20 & 12.975\\
\bottomrule
\end{tabular}}
\end{subtable}
 \par\medskip 
\begin{subtable}{0.48\textwidth}
\centering
\caption{Outlier-contaminated noise}
\adjustbox{max width=\linewidth}{
\begin{tabular}{c|ccc|cc|cc|cc}
\toprule
$c$ & RUKF & RCKF & REKF & $N_p$ & PF & $\tau$ & ME-PF & $\sigma$ & MC-UKF \\
\midrule
0      & 4.1364 & 5.0350 & 12.731 & 100   & 7.8873 & 0.0001 & 10.371      & 1  & 13.357 \\
0.0001 & 3.8280 & 4.1585 & 80.043 & 500   & 7.3860 &  0.01 &9.6539  & 2  & 12.028 \\
0.001  & 3.3129 & 2.7156 & 25.581 & 1000  & 7.1357 & 0.03 & 9.0540 & 4  & 10.947 \\
0.01   & 2.1777 & 1.2692 & 107.77 & 2000  & 6.9964 & 0.05 &  8.8862 & 6  & 10.670 \\
0.03   & 1.7433 & 0.9874 & 2.6127 & 5000  & 6.8166 & 0.1  &8.4875   & 8  & 10.730 \\
0.05   & \textbf{1.6610} & 0.9332 & 2.4449 & 10000 & 6.5092 & 0.3 & \textbf{8.1744}  & 10 & \textbf{10.508} \\
0.1    & 1.8094 & \textbf{0.9177} & \textbf{2.3561} & 20000 & 6.3178 & 0.5 & 8.2900  & 15 & 11.199 \\
0.3    & 2.8316 & 0.93888 & 2.6141 & 50000 &    \textbf{6.1497}    & 0.8  & 8.4104 & 20 & 11.498 \\
\bottomrule
\end{tabular}}
\end{subtable}
\hfill
\begin{subtable}{0.48\textwidth}
\centering
\caption{Dead-zone nonlinearity}
\adjustbox{max width=\linewidth}{
\begin{tabular}{c|ccc|cc|cc|cc}
\toprule
$c$ & RUKF & RCKF & REKF & $N_p$ & PF & $\tau$ & ME-PF & $\sigma$ & MC-UKF \\
\midrule
0      & 3.7425 & 4.4650 & \textbf{4.9081} & 100   & 7.0916 &0.0001 &  9.4790   & 1  & 11.562 \\
0.0001 & 3.7425 & 4.4650 & 27.048 & 500   & 6.5861 & 0.01&  8.7157    & 2  & 10.775 \\
0.001  & 2.9910 & 2.4146 & 91.527 & 1000  & 6.3933 &0.03& 8.3464  & 4  & 9.8940 \\
0.01   & 1.9826 & 1.1419 & 40.957 & 2000  & 6.2419 &0.05&7.8454& 6  & 9.6685 \\
0.03   & 1.6033 & 0.8783 & 26.681 & 5000  & 5.9933 &0.1& 7.5429  & 8  & 9.6023 \\
0.05   & \textbf{1.5416} & 0.8011 & 160.94 & 10000 & 5.8311 &0.3& \textbf{7.3655}    & 10 & \textbf{9.5392} \\
0.1    & 1.7125 & \textbf{0.7554} & 644.02 & 20000 & 5.6987 &0.5&  7.4180 & 15 & 10.598 \\
0.3    & 2.7847 & 0.7934 & 1483.4 & 50000 &    \textbf{5.5102}  & 0.8& 7.6307    & 20 & 11.211 \\
\bottomrule
\end{tabular}}
\end{subtable}
\label{tab1}
\end{table*}
\begin{table}[t]
\centering
\caption{$\overline{\text{MSE}}$ for the RUKF  with $c=0.03$  and UKF for different values of $\alpha$ under the sensor drift scenario.}
\begin{tabular}{|c|c|c|c|c|}
\hline
$\alpha$ & 0.6 & 0.7 & 0.8 & 1 \\ \hline
RUKF & 1.9506 & 1.6790 & 1.5736 & 1.5407 \\ \hline
UKF  & 4.1688 & 4.2860 & 4.5367 & 4.7884 \\ \hline
\end{tabular}
\label{tab_alpha}
\end{table}  we consider the following nonlinear filtering algorithms: the standard EKF, UKF, CKF,  the bootstrap particle filter (PF) with $N_p$ particles~\cite{sarkka2023bayesian}, and the proposed robust sigma point filters, i.e. RUKF  and RCKF. 
Moreover, to further highlight the effectiveness of the proposed approach, we conduct a systematic comparison with several state-of-the-art robust nonlinear filtering methods available in the literature. 
\begin{itemize}
\item REKF: the robust extended Kalman filter \cite{longhini2021learning} with the following changes. The ambiguity set is chosen as in \eqref{ambi_st}. 
In addition, since this method is not directly applicable in the present setting due to the non-differentiability of $f(x)$, we compute the numerical Jacobian matrix   via  a Newton approximation method.
\item MC-UKF: the robust maximum correntropy UKF~\cite{zhao2022robust}, where $\sigma$ denotes the kernel width and the parameters in~(\ref{W_ukf}) are set identically to those of  UKF and RUKF.  The implementation follows the code publicly made available online by the authors.

  \item ME-PF:  the robust maximum entropy particle filter  proposed  in \cite{9872130} with $1000$ particles. 
  For each particle, the corresponding maximum entropy distribution is assumed to lie within a Kullback--Leibler divergence ball centered at the nominal distribution, with radius denoted by~$\tau$.   The implementation follows the code publicly made available online by the authors.
\end{itemize}      All these filters are initialized with $\hat x_0=[3\; 0]^\top$ and $P_0=0.1 I$.

{\em Results.} The performance of these filters is assessed using the  average of the mean squared error over the entire time horizon 
$\overline{\text{MSE}}_t$  defined in (\ref{def_perf}).
Table~\ref{tab:rmse_comparison} summarizes the performance of various filters under different measurement non-idealities  and filter parameter settings.
Note that when \( c = 0 \), RUKF, RCKF, and REKF reduce to the standard UKF, CKF, and EKF, respectively.
It is clear that all the proposed robust sigma point Kalman filters outperform  the standard filters (i.e. UKF, CKF, EKF, and PF), regardless of the specific value of \( c \).
For sigma point filters, when \( c \) is very small, their performance is comparable to that of the standard filters.
As \( c \) increases, the performance improves; however, for large values of \( c \), a slight performance degradation is observed, as the robust filters become overly conservative.  With regard to REKF, its performance is  worse than that of the proposed estimator. It is worth noting that  REKF, RUKF and RCKF are based on the same robust paradigm and therefore adopt the same strategy for coping with uncertainty, but differ in the way they approximate the nonlinear dynamics. This indicates that, in practice, it is important not only to achieve robustness against uncertainty, but also to preserve an accurate approximation of the nonlinear dynamics.  Moreover, unlike the other robust filters, the performance of REKF fluctuates irregularly as $c$ varies across almost all Monte Carlo experiments. This behavior is primarily attributable to a small number of realizations in which the estimated velocity enters a neighborhood of zero. Since $\eta(x)$ exhibits a discontinuity on the hyperplane $s=0$, the first-order Taylor approximation employed by REKF becomes highly inaccurate in these realizations, resulting in catastrophic estimation errors that significantly affect the overall performance metrics. Then, PF fails to achieve a good performance due to its high sensitivity to model uncertainties, even when a large number of particles is employed.
The unsatisfactory performance of ME-PF can be attributed to the fact that its notion of robustness is achieved by maximizing the entropy of the particles ignoring the  potential weakness of the robust estimator.
Finally, the proposed robust filters also outperform MC-UKF. One possible reason is that the latter does not incorporate a time-varying strategy for adjusting the kernel width $\sigma$ at each time step. 
To further investigate the impact of the parameter $\alpha$ on the performance of RUKF, different values of $\alpha$ are selected under the sensor drift scenario, as summarized in Table \ref{tab_alpha}: changing parameter setting does not alter the fact that the RUKF performs better than the UKF.  Finally, we compare the computational time of the robust filters. 
All simulations were implemented in MATLAB and executed on a Mechanical Revolution notebook 
equipped with an AMD R9-7845HX CPU and a GeForce RTX 4070 GPU. 
To ensure a fair comparison, fixed parameter values (i.e., $c$, $\tau$, or $\sigma$) 
were used for all the robust filters, chosen as the mean of the candidate values listed in 
Table~\ref{tab1}.  As shown in Table \ref{runningt},  the proposed robust sigma-point Kalman filters have  running time similar to MC-UKF and REKF, 
yet are significantly more computationally efficient than ME-PF.


\begin{table}[t]
\centering
\caption{Mean running time of robust filters across four different scenarios.}
\label{runningt}
\small

\resizebox{\columnwidth}{!}{%
\begin{tabular}{c c c c c}
\hline
Method & Sensor drift & Intermittent obs. & Outlier cont. noise & Nonlinearity \\
\hline
RUKF  & 0.0030\,s & 0.0029\,s & 0.0026\,s & 0.0026\,s \\
RCKF  & 0.0029\,s & 0.0029\,s & 0.0026\,s &  0.0027\,s \\
REKF  & 0.0021\,s & 0.0019\,s &0.0017 \,s & 0.0017\,s \\
ME-PF  & 1.5915\,s & 1.6287\,s & 1.6372\,s & 1.5843\,s \\
MC-UKF  & 0.0029\,s &  0.0028\,s & 0.0029\,s &0.0027\,s \\
\hline
\end{tabular}
}
\end{table}

\section{Conclusion} \label{sec_conc}
We proposed a robust sigma-point filter and interpreted it as the minimizer of a minimax game whose ambiguity set is centered on a sigma-point approximation.  We have shown that such an interpretation does not always hold, but it depends on whether function $h$ satisfies certain properties.
Moreover, we  proposed a MCMC scheme for generating adversarial data from the least favorable model, thus  allowing uncertainty to be assessed by examining the resulting realizations. The numerical results confirmed the convergence of the proposed MCMC scheme and validated the expected behavior of the robust sigma-point filter under worst-case conditions. Finally, they also showed the effectiveness of the proposed filtering approach compared to the existing robust nonlinear filters.


\begin{thebibliography}{10}

\bibitem{reif1999stochastic}
K.~Reif, S.~Gunther, E.~Yaz, and R.~Unbehauen, ``Stochastic stability of the
  discrete-time extended {K}alman filter,'' {\em IEEE Transactions on Automatic
  Control}, vol.~44, no.~4, pp.~714--728, 1999.

\bibitem{julier2000new}
S.~Julier, J.~Uhlmann, and H.~F. Durrant-Whyte, ``A new method for the
  nonlinear transformation of means and covariances in filters and
  estimators,'' {\em IEEE Transactions on Automatic Control}, vol.~45, no.~3,
  pp.~477--482, 2000.

\bibitem{arasaratnam2009cubature}
I.~Arasaratnam and S.~Haykin, ``Cubature {K}alman filters,'' {\em IEEE
  Transactions on Automatic Control}, vol.~54, no.~6, pp.~1254--1269, 2009.

\bibitem{pruher2020improved}
J.~Pr{\"u}her, T.~Karvonen, C.~J. Oates, O.~Straka, and S.~S{\"a}rkk{\"a},
  ``Improved calibration of numerical integration error in sigma-point
  filters,'' {\em IEEE Transactions on Automatic Control}, vol.~66, no.~3,
  pp.~1286--1292, 2020.

\bibitem{arasaratnam2007discrete}
I.~Arasaratnam, S.~Haykin, and R.~J. Elliott, ``Discrete-time nonlinear
  filtering algorithms using {G}auss--{H}ermite quadrature,'' {\em Proceedings
  of the IEEE}, vol.~95, no.~5, pp.~953--977, 2007.

\bibitem{sarkka2023bayesian}
S.~S{\"a}rkk{\"a} and L.~Svensson, {\em Bayesian filtering and smoothing},
  vol.~17.
\newblock Cambridge university press, 2023.

\bibitem{cheng2026distributionally1}
J.~Cheng, Z.~Xue, H.~Chen, and Y.~Huang, ``Distributionally robust {K}alman
  filter under likelihood model mismatch from an optimization perspective-part
  ii: Extension and application to cooperative localization,'' {\em IEEE
  Transactions on Aerospace and Electronic Systems}, 2026.

\bibitem{cheng2026distributionally2}
J.~Cheng, Z.~Xue, H.~Chen, and Y.~Huang, ``Distributionally robust {K}alman
  filter under likelihood model mismatch from an optimization perspective-part
  ii: Extension and application to cooperative localization,'' {\em IEEE
  Transactions on Aerospace and Electronic Systems}, 2026.

\bibitem{8025799}
Y.~Huang, Y.~Zhang, Z.~Wu, N.~Li, and J.~Chambers, ``A novel adaptive {K}alman
  filter with inaccurate process and measurement noise covariance matrices,''
  {\em IEEE Transactions on Automatic Control}, vol.~63, no.~2, pp.~594--601,
  2018.

\bibitem{9146725}
Y.~Huang, Y.~Zhang, Y.~Zhao, P.~Shi, and J.~A. Chambers, ``A novel
  outlier-robust {K}alman filtering framework based on statistical similarity
  measure,'' {\em IEEE Transactions on Automatic Control}, vol.~66, no.~6,
  pp.~2677--2692, 2021.

\bibitem{nakabayashi2019nonlinear}
A.~Nakabayashi and G.~Ueno, ``Nonlinear filtering method using a switching
  error model for outlier-contaminated observations,'' {\em IEEE Transactions
  on Automatic Control}, vol.~65, no.~7, pp.~3150--3156, 2019.

\bibitem{zhao2022robust}
H.~Zhao, B.~Tian, and B.~Chen, ``Robust stable iterated unscented {K}alman
  filter based on maximum correntropy criterion,'' {\em Automatica}, vol.~142,
  p.~110410, 2022.

\bibitem{ROBUST_STATE_SPACE_LEVY_NIKOUKHAH_2013}
B.~Levy and R.~Nikoukhah, ``Robust state-space filtering under incremental
  model perturbations subject to a relative entropy tolerance,'' {\em {IEEE
  Transactions on Automatic Control}}, vol.~58, pp.~682--695, Mar. 2013.

\bibitem{STATETAU_2017}
M.~Zorzi, ``Robust {K}alman filtering under model perturbations,'' {\em {IEEE
  Transactions on Automatic Control}}, vol.~62, no.~{6}, pp.~{2902--2907},
  2016.

\bibitem{kim2020robust}
S.~Kim, V.~Deshpande, and R.~Bhattacharya, ``Robust {K}alman filtering with
  probabilistic uncertainty in system parameters,'' {\em {IEEE Control Systems
  Letters}}, vol.~5, no.~1, pp.~295--300, 2020.

\bibitem{abadeh2018wasserstein}
S.~Abadeh, V.~Nguyen, D.~Kuhn, and P.~Esfahani, ``Wasserstein distributionally
  robust kalman filtering,'' in {\em Advances in Neural Information Processing
  Systems}, pp.~8474--8483, 2018.

\bibitem{robustleastsquaresestimation}
B.~Levy and R.~Nikoukhah, ``Robust least-squares estimation with a relative
  entropy constraint,'' {\em {IEEE Transactions on Information Theory}},
  vol.~50, no.~1, pp.~89--104, 2004.

\bibitem{yi2021robust}
S.~Yi and M.~Zorzi, ``Robust {K}alman filtering under model uncertainty: The
  case of degenerate densities,'' {\em IEEE Transactions on Automatic Control},
  vol.~67, no.~7, pp.~3458--3471, 2021.

\bibitem{10654520}
Y.~Xu, W.~Xue, C.~Shang, and H.~Fang, ``On globalized robust {K}alman filter
  under model uncertainty,'' {\em {IEEE Transactions on Automatic Control}},
  vol.~70, no.~2, pp.~1147--1160, 2025.

\bibitem{longhini2021learning}
A.~Longhini, M.~Perbellini, S.~Gottardi, S.~Yi, H.~Liu, and M.~Zorzi,
  ``Learning the tuned liquid damper dynamics by means of a robust {EKF},'' in
  {\em 2021 American Control Conference (ACC)}, pp.~60--65, IEEE, 2021.

\bibitem{11107768}
S.~Yi, X.~Jin, Z.~Wang, Z.~Liu, and M.~Zorzi, ``Data-driven robust uav position
  estimation in gps signal-challenged environment,'' in {\em 2025 American
  Control Conference (ACC)}, pp.~491--496, 2025.

\bibitem{jang2026residual}
M.~Jang, J.~Lee, A.~Hakobyan, N.~Hovakimyan, and I.~Yang, ``Residual-aware
  distributionally robust {EKF}: Absorbing linearization mismatch via
  wasserstein ambiguity,'' {\em arXiv preprint arXiv:2604.02749}, 2026.

\bibitem{URKF}
S.~Yi and M.~Zorzi, ``An update-resilient {K}alman filtering approach,'' {\em
  Accepted in Automatica}, 2026.

\bibitem{9872130}
S.~Wang, ``Distributionally robust state estimation for nonlinear systems,''
  {\em IEEE Transactions on Signal Processing}, vol.~70, pp.~4408--4423, 2022.

\bibitem{yi2026eccspKF}
S.~Yi and M.~Zorzi, ``Robust sigma-point filtering for nonlinear systems with
  non-additive noise,'' in {\em 2026 European Control Conference (ECC),
  accepted}, IEEE, 2026.

\bibitem{hansen2008robustness}
L.~P. Hansen and T.~J. Sargent, {\em Robustness}.
\newblock Princeton university press, 2008.

\bibitem{aubin2006applied}
J.~Aubin and I.~Ekeland, {\em Applied nonlinear analysis}.
\newblock Courier Corporation, 2006.

\bibitem{wan2001unscented}
E.~A. Wan and R.~Van Der~Merwe, ``The unscented kalman filter,'' {\em Chapter 7
  of: Kalman filtering and neural networks}, pp.~221--280, 2001.

\bibitem{zenere2018coupling}
A.~Zenere and M.~Zorzi, ``On the coupling of model predictive control and
  robust {K}alman filtering,'' {\em IET Control Theory \& Applications},
  vol.~12, no.~13, pp.~1873--1881, 2018.

\bibitem{zorzi2017convergence}
M.~Zorzi, ``Convergence analysis of a family of robust kalman filters based on
  the contraction principle,'' {\em SIAM Journal on Control and Optimization},
  vol.~55, no.~5, pp.~3116--3131, 2017.

\bibitem{hill2019stationarity}
S.~D. Hill and J.~C. Spall, ``Stationarity and convergence of the
  {M}etropolis-{H}astings algorithm: Insights into theoretical aspects,'' {\em
  IEEE Control Systems Magazine}, vol.~39, no.~1, pp.~56--67, 2019.

\bibitem{khalil2002nonlinear}
H.~K. Khalil and J.~W. Grizzle, {\em Nonlinear systems}, vol.~3.
\newblock Prentice hall Upper Saddle River, NJ, 2002.

\end{thebibliography}
\end{document}